\documentclass[12pt,reqno]{amsart}
\usepackage{amssymb}
\usepackage{amsmath}
\usepackage{amsfonts}

\theoremstyle{plain}

\numberwithin{equation}{section}

\begin{document}
\title[On a non local non-homogeneous fractional Timoshenko system]{On a non
local non-homogeneous fractional Timoshenko system with frictional and
viscoelastic damping terms}
\author{Said Mesloub, Eman Alhazzani, Gadain Hassan Eltayeb }
\address{Mathematics Department, College of Science, King Saud University,
P.O. Box 2455, Riyadh 11451, Saudi Arabia}
\email{mesloub@ksu.edu.sa; \ gadain@ksu.edu.sa; \ esalhazani@imamu.edu.sa}

\begin{abstract}
We are devoted to the study of a nonhomogeneous time fractional timoshenko
system with frictional and viscoelastic damping terms. We are concerned with
the well-posedness of the given problem. The approach relies on some
fucntional analysis tools, operator theory, a prori estimates and density
arguments.
\end{abstract}

\maketitle

\section{Introduction}

Vibrations of beams are not always safe and welcomed because of their great
and irreparable damages effects. In this situation, researchers try to
introduce some damping mechanisms (viscous damping, thermoelastic damping,
modal damping, frictional damping) in such a way that these damaging and
destructive vibrations are perfectly reduced. In other words, an intensive
investigation has been carried out to impose minimal conditions to provide
and guarantee stability of Timoshenko systems using several types of
dissipative mechanisms.

As a classical and a simple model [1], Timoshenko studied the following
coupled hyperbolic system%
\begin{equation}
\left\{
\begin{array}{c}
\begin{array}{c}
\rho _{1}\theta _{tt}-\kappa (\theta _{x}-\mathcal{\phi })_{x}=0,\text{ \ \
\ \ \ \ \ \ \ \ \ \ \ \ }(x,t)\in \text{ }(0,L)\times (0,\infty ) \\
\rho _{2}\mathcal{\phi }_{tt}=\kappa ^{\ast }\mathcal{\phi }_{xx}+\kappa
(\theta _{x}-\mathcal{\phi })\text{ \ \ \ }(x,t)\in \text{ }(0,L)\times
(0,\infty ),%
\end{array}
\\
(\theta _{x}-\mathcal{\phi })\mid _{x=0}^{x=L}=0,\text{ \ }\mathcal{\phi }%
_{x}\mid _{x=0}^{x=L}=0,%
\end{array}%
\right.  \label{e1.1*}
\end{equation}%
describing the transverse vibration of a beam. where $L$ is the length of
the beam in its equilibrium configuration. The function $\theta $ models the
transverse displacement of the beam and $\mathcal{\phi }$ models the
rotation angle of its filament. The coefficients $\rho _{1},$ $\rho _{2},$ $%
\kappa $ and $\kappa ^{\ast }$ are respectively the density, the polar
moment of inertia of a cross section, the shear modulus and \ the Young's
modulus of elasticity. Timoshenko system (\ref{e1.1*}) was generalized and
studied by many authors. As mentioned at the beginning of the introduction,
different types of dampings were added to the Timoshenko system for the
purpose of its stabilization. For example, in [4], researchers investigated
the exponential stability for a Timoshenko \ system having two weak dampings%
\begin{equation}
\left\{
\begin{array}{c}
\rho _{1}\theta _{tt}=\kappa (\theta _{x}-\mathcal{\phi })_{x}-\theta _{t},%
\text{ \ \ in }(0,L)\times (0,\infty ), \\
\rho _{2}\mathcal{\phi }_{tt}=\kappa ^{\ast }\mathcal{\phi }_{xx}-\kappa
(\theta _{x}-\mathcal{\phi })_{x}-\mathcal{\phi }_{t},\text{ \ \ in }%
(0,L)\times (0,\infty ), \\
\theta (0,t)=\theta (L,t)=\mathcal{\phi }(0,t)=\mathcal{\phi }(L,t)=0,\text{
\ }t>0.%
\end{array}%
\right.  \label{e1.2*}
\end{equation}

In [2], authors proved some exponential decay results for a Timoshenko
system with a memory damping term%
\begin{equation}
\left\{
\begin{array}{c}
\rho _{1}\theta _{tt}-\kappa _{1}(\theta _{x}+\mathcal{\phi })_{x}=0,\text{
\ \ \ in }(0,L)\times (0,\infty ) \\
\rho _{2}\mathcal{\phi }_{tt}-\kappa _{2}\mathcal{\phi }_{xx}+\kappa
_{1}(\theta _{x}+\mathcal{\phi })+h\ast \mathcal{\phi }_{xx}(x,t)=0,\text{
in }(0,L)\times (0,\infty ) \\
\theta (0,t)=\theta (L,t)=\mathcal{\phi }(0,t)=\mathcal{\phi }(L,t)=0, \\
\theta (x,0)=\theta _{0},\text{ }\theta _{t}(x,0)=\theta _{1},\text{ }%
\mathcal{\phi }(x,0)=\mathcal{\phi }_{0},\text{ }\mathcal{\phi }_{t}(x,0)=%
\mathcal{\phi }_{1}.%
\end{array}%
\right.  \label{e1.3*}
\end{equation}%
Authors considered and studied in [5] the effect of frictional and
viscoelastic dampings, and proved some exponential and polynomial decay
results for the system%
\begin{equation}
\left\{
\begin{array}{c}
\theta _{tt}-(\theta _{x}+\mathcal{\phi })_{x}=0, \\
\mathcal{\phi }_{tt}-\mathcal{\phi }_{xx}+\theta _{x}+\mathcal{\phi }%
+\int\limits_{0}^{t}g(t-s)(a(x)\mathcal{\phi }_{x}(x,s))_{x}ds+b(x)h(%
\mathcal{\phi }_{t})=0, \\
\theta (0,t)=\theta (1,t)=\mathcal{\phi }(0,t)=\mathcal{\phi }(1,t)=0,\text{
\ }t>0.%
\end{array}%
\right.  \label{e1.4*}
\end{equation}%
We also mention that in [$44$], the authors investigated the exponential
stabilization of a Timoshenko system by a thermal effect damping.%
\begin{equation}
\left\{
\begin{array}{c}
\rho _{1}\theta _{tt}-\kappa _{1}(\theta _{x}+\mathcal{\phi })_{x}=0,\text{
\ \ \ in }(0,L)\times (0,\infty ) \\
\rho _{2}\mathcal{\phi }_{tt}-\kappa _{2}\mathcal{\phi }_{xx}+\kappa
_{1}(\theta _{x}+\mathcal{\phi })+\gamma \mathcal{\omega }_{x},\text{ in }%
(0,L)\times (0,\infty ) \\
\rho _{2}\mathcal{\omega }_{tt}-\kappa _{3}\mathcal{\omega }_{xx}+\beta
\int\limits_{0}^{t}g(t-s)\mathcal{\omega }_{xx}(x,s)ds+\gamma \mathcal{\phi }%
_{ttx},\text{ in }(0,L)\times (0,\infty ).%
\end{array}%
\right.  \label{e1.5*}
\end{equation}%
In [3], the author considered a \ Timoshenko linear thermoelastic system
with linear frictional damping and a distributed delay. He proved the
well-posedness, and proved that the system is exponentially stable
regardless of the speeds of wave propagation.There are many other papers in
the literature dealing with the stabilization of different version of
Timoshenko systems. For more results concerning the stabilization and
controllability of Timoshenko systems, we refer the reader to [6, 7, 8, 9),
10, 11, 12, 13, 14, 15, 31].

Recently, a generalization of the Timoshenko system (\ref{e1.1*}) into
fractional setting is studied in [45] by using a fractional version of
resolvents. The author established the well posedness of a fractional
Timoshenko system, and proved that lower order fractional terms can
stabilize the system in a Mittag-Leffler fashion. More precisely, the author
considered the initial boundary value problem%
\begin{equation*}
\left\{
\begin{array}{c}
\rho _{1}\partial _{t}^{\alpha }(\partial _{t}^{\alpha }\theta )-\kappa
_{1}(\theta _{x}+\mathcal{\phi })_{x}=0,\text{ \ \ \ in }(0,1)\times
(0,\infty ) \\
\rho _{2}\partial _{t}^{\alpha }(\partial _{t}^{\alpha }\mathcal{\phi +}a%
\mathcal{\phi })-\kappa _{2}\mathcal{\phi }_{xx}+\kappa _{1}(\theta _{x}+%
\mathcal{\phi }),\text{ in }(0,1)\times (0,\infty ) \\
\theta (0,t)=\theta (1,t)=0,\text{ }\mathcal{\phi }(0,t)=\mathcal{\phi }%
(1,t)=0,\text{ \ }t>0 \\
\theta (x,0)=\theta _{0}(x)\ ,\ \mathcal{\phi (}x,0)=\psi (x).%
\end{array}%
\right.
\end{equation*}

Motivated by the above results on Timoshenko systems, we consider a non
local initial boundary value problem for a non-homogeneous fractional
Timoshenko system with a frictional damping in the first equation and a
viscoelastic memory damping term in the second equation. The system is
complemented with initial conditions and non local purely boundary integral
conditions. At the beginning of the year 1963, Cannon [20] was the first
researcher to investigate a non local problem with a non local constraint
(energy specification) of the form $\int\limits_{0}^{l}\chi (x)U(x,t)dt=\tau
(t)$, where $\chi (x),$ and $\tau (t)$ are given functions. More precisely,
he used the potential method to investigate the well posedness of the heat
equation subject to the specification of energy. This type of conditions
arise mainly when the data cannot be measured directly on the boundary, but
only their averages (weighted averages) are known. Due to their importance,
physical significance (mean, total flux, total energy,..) and numerous
applications in different fields of science and engineering, such as
underground water flow, vibration problems, heat conduction, medical
science, nuclear reactor dynamics, thermoelasticity, and plasma physics and
control theory, several authors extensively studied this type of problems.
We can cite for example [16, 21, 22, 23, 24, 25, 26, $28,29,30,31$]. Note
that theoretical study of non local problems is connected with great
difficulties, since the presence of integral terms in the boundary
conditions can greatly complicate the application of classical methods of
functional analysis method, especially when it comes to the fractional case
. A functional analysis method based on some a priori bounds and on the
density of the range of the unbounded operator corresponding to the abstract
formulation of the given problem is used to prove the well posedness of the
posed problem. This is shown through the introduction of some multiplier
operators, some classical and fractional inequalities, and the establishment
of some properties, involving fractional derivatives.

To the best of our knowledge, the treated fractional system problem (\ref%
{e2.1})-(\ref{e2.4}) has never been studied and explored in the literature.
This work can be considered as a contribution in the development of the
traditional functional analysis method, the so called energy inequality
method used to prove the well posedness of mixed problems with integral
boundary conditions. For some classical cases, the reader can refeer for
example to example [16, 17, 18, 19, 27], and for some fractional cases, the
reader should refeer to $[32,33,34,36,37,38,39,40,41].$ We should also
mention here that there are some important papers dealing with numerical
aspects for Timoshenko systems, and having many applications, for which the
reader can refer to [47 48, 49, 50]. There are some papers dealing with
Timoshenko system with fractional operator in the memory [42, 43].

\section{Formulation of the problem and function spaces}

Given the interval $I=(0,L),$ we consider the non-homogeneous fractional
viscoelastic beam model with frictional damping of Timoshenko type%
\begin{equation}
\left\{
\begin{array}{c}
\mathcal{L}_{1}(\theta ,\mathcal{\phi })=\rho _{1}\partial _{t}^{\alpha
+1}\theta -\kappa _{1}(\theta _{x}+\mathcal{\phi })_{x}+\theta _{t}=F(x,t)
\\
\mathcal{L}_{2}(\theta ,\mathcal{\phi })=\rho _{2}\partial _{t}^{\alpha +1}%
\mathcal{\phi }-\kappa _{2}\mathcal{\phi }_{xx}+\kappa _{1}(\theta _{x}+%
\mathcal{\phi })+\int\limits_{0}^{t}m(t-s)\mathcal{\phi }_{xx}(x,s)ds=G(x,t),%
\end{array}%
\right.  \label{e2.1}
\end{equation}%
in the unknowns $(\theta ,\phi ):(x,t)$ $\in I\times \lbrack 0,T]\rightarrow
\mathbb{R}
,$ the strictly positive constants $\rho _{1},\rho _{2},\kappa _{1}$ and $%
\kappa _{2}$ satisfy the relation%
\begin{equation*}
\frac{\rho _{1}}{\kappa _{1}}=\frac{\rho _{2}}{\kappa _{2}},
\end{equation*}%
and $\ f\ ,\ g\ ,\ \varphi \ ,\ \psi ,$ $F,$ and $G\ $are given functions,
and $m:$ $%
\mathbb{R}
^{+}\rightarrow
\mathbb{R}
^{+}$ is a twice differentiable function such that%
\begin{equation}
\kappa _{2}-\int\limits_{0}^{T}m(t)dt=l>0,\text{ \ }m^{\prime }(t)<0,\text{ }%
\forall t\geq 0.  \label{e2.2}
\end{equation}%
The system (\ref{e2.1}) is complemented with the initial conditions%
\begin{equation}
\left\{
\begin{array}{c}
\Gamma _{1}\theta =\theta (x,0)=\varphi (x)\ ,\ \Gamma _{2}\theta =\theta
_{t}(x,0)=\psi (x)\ , \\
\Gamma _{1}\mathcal{\phi }=\mathcal{\phi }(x,0)=f(x)\ ,\ \Gamma _{2}\mathcal{%
\phi }=\mathcal{\phi }_{t}(x,0)=g(x),%
\end{array}%
\right.  \label{e2.3}
\end{equation}%
and the non local boundary integral conditions%
\begin{equation}
\int\limits_{0}^{L}\theta dx=0\ ,\ \int\limits_{0}^{L}x\theta
dx=0,\int\limits_{0}^{L}\mathcal{\phi }dx=0\ ,\ \int\limits_{0}^{L}x\mathcal{%
\phi }dx=0.  \label{e2.4}
\end{equation}%
This system of coupled hyperbolic equations represents a Timoshenko model
for a thick beam of length $L,$where $\theta $ is the transverse
displacement of the beam and $\mathcal{\phi }$ is the rotation angle of the
filament of the beam. The coefficients $\rho _{1},\rho _{2},\kappa _{1}$ and
$\kappa _{2}$ are respectively the density, the polar moment of inertia of a
cross section, the shear modulus and \ the Young's modulus of elasticity.
The integral conditions represent the averages (weighted averages) of the
total transverse displacement of the beam and the rotation angle of the
filament of the beam.

Our aim is to study the well posedness of the solution of problem (\ref{e2.1}%
), (\ref{e2.4}). That is on the basis of some a priori bounds and on the
density of the range of the operator generated by the problem under
consideration, we prove the existence, uniqueness and continuous dependence
of the solution on the given data of problem (\ref{e2.1}), (\ref{e2.4}). We
now introduce some function spaces needed throughout the sequel. Let $%
L^{2}(Q^{T})$ be the Hilbert space of square integrable functions on $%
Q^{T}=(0,1)\times (0,T),$ $T<\infty ,$ with scalar product and norm
respectively%
\begin{equation}
(Z,S)_{L^{2}(Q^{T})}=\int_{Q^{T}}ZSdxdt,\text{ \ \ \ }\left\Vert
Z\right\Vert _{L^{2}(Q^{T})}^{2}=\int_{Q^{T}}Z^{2}dxdt.  \label{e2.5}
\end{equation}%
We also use the space $L^{2}((0,1))$ on the interval $(0,1)$, whose
definition is analogous to the space on $Q.$ Let $B_{2}^{1}(0,L)$ be the
space obtained by completion of the space $C^{0}(0,L)$ of real continuous
functions with compact support in the interval $(0,L)$ with respect to the
inner product%
\begin{equation*}
(\gamma ,\gamma ^{\ast })_{B_{2}^{1}(0,L)}=\int\limits_{0}^{L}\mathcal{I}%
_{x}\gamma .\Im _{x}\gamma ^{\ast }dx,
\end{equation*}%
where $\mathcal{I}_{x}\gamma =\int\limits_{0}^{x}\gamma (\zeta )d\zeta $ for
every fixed $x\in (0,L).$ The associated norm is $\left\Vert \gamma
\right\Vert _{B_{2}^{1}(0,L)}^{2}=\sqrt{(\gamma ,\gamma )_{B_{2}^{1}(0,L)}}%
=\int\limits_{0}^{L}\left( \mathcal{I}_{x}\gamma \right) ^{2}dx.$ We denote
by $C(\overline{J};L^{2}(0,L))$ with $J=(0,T)$ the set of all continuous
functions $\gamma (.,t):J\rightarrow L^{2}(0,L)$ with norm%
\begin{equation}
\left\Vert \gamma \right\Vert _{C(J;L^{2}(0,L))}^{2}=\underset{0\leq t\leq T}%
{\sup }\left\Vert \gamma (.,t)\right\Vert _{L^{2}(0,L)}^{2}<\infty ,
\label{e2.6}
\end{equation}%
and $C(\overline{J};B_{2}^{1}(0,L))$ the set of functions $\gamma (.,t):%
\overline{J}\rightarrow B_{2}^{1}(0,L)$ with norm%
\begin{equation}
\left\Vert \gamma \right\Vert _{C(\overline{J};B_{2}^{1}(0,L))}^{2}=\underset%
{0\leq t\leq T}{\sup }\left\Vert \mathcal{I}_{x}\gamma (.,t)\right\Vert
_{L^{2}(0,L)}^{2}=\underset{0\leq t\leq T}{\sup }\left\Vert \gamma
(.,t)\right\Vert _{B_{2}^{1}(0,L)}^{2}<\infty .  \label{e2.7}
\end{equation}%
To obtain \ a priori estimate for the solution, we write down our problem (%
\ref{e2.1}), (\ref{e2.4}) in its operator form: $\mathcal{GZ}=H$ with $%
\mathcal{Z}=(\mathcal{\theta },\phi )$, $\mathcal{GZ}=(\mathcal{S}_{1}(%
\mathcal{\theta },\phi ),\mathcal{S}_{2}(\mathcal{\theta },\phi ))$ and $%
H=(H_{1},H_{2})$ where%
\begin{equation}
\left\{
\begin{array}{c}
\begin{array}{c}
\mathcal{S}_{1}(\mathcal{\theta },\phi )=\{\mathcal{L}_{1}(\mathcal{\theta }%
,\phi ),\Gamma _{1}\mathcal{\theta },\Gamma _{2}\mathcal{\theta }\} \\
\mathcal{S}_{2}(\mathcal{\theta },\phi )=\{\mathcal{L}_{2}(\mathcal{\theta }%
,\phi ),\Gamma _{1}\phi ,\Gamma _{2}\phi \}%
\end{array}
\\
H_{1}=\{F,\varphi ,\psi \},\text{ }H_{2}=\{G,f,g\}.%
\end{array}%
\right. ,\text{ }  \label{e2.8}
\end{equation}%
The operator $\mathcal{G}$ is an unbounded operator of domain of definition $%
D(\mathcal{G})$ consisting of elements $(\mathcal{\theta },\phi )\in \left(
L^{2}(\overline{J};L^{2}(0,L))\right) ^{2}$ such that $\mathcal{\theta }%
_{x},\phi _{x},\mathcal{\theta }_{t},\phi _{t},\mathcal{\theta }_{tt},\phi
_{tt},\mathcal{\theta }_{xx},\phi _{xx}$ belonging to $L^{2}(\overline{J}%
;L^{2}(0,L))$ verifying initial and boundary conditions (\ref{e2.3}) and (%
\ref{e2.4}). The operator $\mathcal{G}$ is acting from the Banach space $%
\mathcal{B}$ into the Hilbert space $\mathcal{E}$, where $\mathcal{B}$ is
the Banach space obtained by completing $D(\mathcal{G})$ with respect to the
norm%
\begin{equation}
\Vert \mathcal{Z}\Vert _{\mathcal{B}}^{2}=\Vert \mathcal{\theta }(.,t)\Vert
_{C(\overline{J};L^{2}(0,L))}^{2}+\Vert \phi (.,t)\Vert _{C(\overline{J}%
;L^{2}(0,L))}^{2}.  \label{e2.9*}
\end{equation}%
And $\mathcal{E}=\left[ L^{2}(Q^{T})\times (L^{2}(0,L))^{2}\right] \times %
\left[ L^{2}(Q^{T})\times (L^{2}(0,L))^{2}\right] $ is the Hilbert space
consisting of vector-valued functions $H=\left( \{F,\varphi ,\psi
\},\{G,f,g\}\right) $ for which the norm
\begin{eqnarray}
\Vert H\Vert _{\mathcal{E}}^{2} &=&\Vert F\Vert _{L^{2}(Q^{T})}^{2}+\Vert
\varphi \Vert _{L^{2}(0,L)}^{2}+\Vert \psi \Vert _{L^{2}(0,L)}^{2}+\Vert
G\Vert _{L^{2}(Q^{T})}^{2}  \notag \\
&&+\Vert f\Vert _{L^{2}(0,L)}^{2}+\Vert g\Vert _{L^{2}(0,L)}^{2}.
\label{e2.10}
\end{eqnarray}%
is finite.

\section{Preliminaries (Definitis and lemmas)}

In this section, we provide some definitions and lemmas needed for
establishing different proves in the sequel.

\textbf{Definition 1.} [51] The time fractional derivative of order $\beta ,$
with $\beta \in (1,2)$ for a function $V$ is defined by%
\begin{equation}
^{C}\partial _{t}^{\beta }V(x,t)=\frac{1}{\Gamma (2-\beta )}%
\int\limits_{0}^{t}\frac{V_{\tau \tau }(x,\tau )}{(t-\tau )^{\beta -1}}d\tau
,  \label{e2.11*}
\end{equation}%
and for $\beta \in (0,1)$ it is defined by%
\begin{equation*}
^{C}\partial _{t}^{\beta }V(x,t)=\frac{1}{\Gamma (1-\beta )}%
\int\limits_{0}^{t}\frac{V_{\tau }(x,\tau )}{(t-\tau )^{\beta }}d\tau
\end{equation*}%
where $\Gamma (1-\beta )$ is the Gamma function.

\textbf{Definition 2.} [51]. The fractional Riemann-Liouville integral of
order $0<\beta <1$ [$52$] which is given by%
\begin{equation}
D_{t}^{-\beta }\upsilon (x,t)=\frac{1}{\Gamma (\beta )}\int\limits_{0}^{t}%
\frac{\upsilon (x,\tau )}{(t-\tau )^{1-\beta }}d\tau .  \tag{2.6}
\end{equation}

\textbf{Lemma 2.1 }$[35]$. Let $E(s)$ be nonnegative and absolutely
continuous on $[0,T]$, and suppose that for almost all $s\in \lbrack 0,T]$, $%
R$ satisfies the inequality%
\begin{equation}
\frac{dE}{ds}\leq A(s)E(s)+B(s),
\end{equation}%
where the functions $A(s)$ and $B(s)$ are summable and nonnegative on $%
[0,T]. $ Then%
\begin{equation}
E(s)\leq \exp \left\{ \int\limits_{0}^{s}A(t)dt\right\} \left(
E(0)+\int\limits_{0}^{s}B(t)dt\right) .
\end{equation}%
\textbf{Lemma 2.3}. $[34]$ Let a nonnegative absolutely continuous function $%
\mathcal{Q}(t)$ satisfy the inequality%
\begin{equation}
^{C}\partial _{t}^{\beta }\mathcal{Q}(t)\leq b_{1}\mathcal{Q}(t)+b_{2}(t),%
\text{ \ \ \ }0<\beta <1,  \label{e2.13*}
\end{equation}%
for almost all $t\in \lbrack 0,T],$ where $b_{1}$ is a positive constant and
$b_{2}(t)$ is an integrable nonnegative function on $[0,T].$ Then%
\begin{equation}
\mathcal{Q}(t)\leq \mathcal{Q}(0)E_{\beta }(b_{1}t^{\beta })+\Gamma (\beta
)E_{\beta ,\beta }(b_{1}t^{\beta })D_{t}^{-\beta }b_{2}(t),  \label{e2.14*}
\end{equation}%
where%
\begin{equation*}
E_{\beta }(x)=\sum_{n=0}^{\infty }\frac{x^{n}}{\Gamma (\beta n+1)}\text{ and
}E_{\beta ,\mu }(x)=\sum_{n=0}^{\infty }\frac{x^{n}}{\Gamma (\beta n+\mu )},
\end{equation*}%
are the MIttag-Leffler functions

\section{A priori estimate and its Consequences}

In this section, we establish an energy inequality from which we deduce the
uniqueness and continuous dependence of solution of problem (\ref{e2.1})-(%
\ref{e2.4}) on the given data.

\textbf{Theorem 3.1}. For any function $\mathcal{Z}=(\mathcal{\theta },\phi
)\in D(\mathcal{G})$ the following a priori estimates holds%
\begin{eqnarray}
&&\Vert \mathcal{\theta }(.,t)\Vert _{C(\overline{J};L^{2}(0,L))}^{2}+\Vert
\phi (.,t)\Vert _{C(\overline{J};L^{2}(0,L))}^{2}  \notag \\
&\leq &\mathcal{F}^{\ast }\left( \Vert \varphi \Vert _{L^{2}(0,L)}^{2}+\Vert
\psi \Vert _{L^{2}(0,L)}^{2}+\Vert g\Vert _{L^{2}(0,L)}^{2}+\Vert f\Vert
_{L^{2}(0,L)}^{2}\right.  \notag \\
&&\left. +\Vert F\Vert _{L^{2}(0,t;L^{2}(0,L))}^{2}+\Vert G\Vert
_{L^{2}(0,t;L^{2}(0,L))}^{2}\right) ,  \label{e3.1}
\end{eqnarray}%
and%
\begin{eqnarray}
&&D^{\alpha -1}\left( \left\Vert \mathcal{\theta }_{t}\right\Vert
_{B_{2}^{1}(0,L)}+\left\Vert \phi _{t}\right\Vert _{B_{2}^{1}(0,L)}\right)
\notag \\
&\leq &\mathcal{F}^{\ast }\left( \Vert \varphi \Vert _{L^{2}(0,L)}^{2}+\Vert
\psi \Vert _{L^{2}(0,L)}^{2}+\Vert g\Vert _{L^{2}(0,L)}^{2}+\Vert f\Vert
_{L^{2}(0,L)}^{2}\right.  \notag \\
&&\left. +\Vert F\Vert _{L^{2}(0,t;L^{2}(0,L))}^{2}+\Vert G\Vert
_{L^{2}(0,t;L^{2}(0,L))}^{2}\right) .  \label{e3.1**}
\end{eqnarray}%
where $\mathcal{F}^{\ast }$ is a positive constant independent of $\mathcal{%
Z=}(\mathcal{\theta },\phi )$ given by%
\begin{equation*}
\mathcal{F}^{\ast }=\mathcal{M}\omega \max \left\{ 1,\frac{T^{\alpha }}{%
\alpha \Gamma (\alpha )}\right\} ,
\end{equation*}%
with%
\begin{eqnarray*}
\mathcal{M} &=&\Gamma (\alpha )E_{\alpha ,\alpha }(\omega t^{\alpha })\left(
\max \left\{ 1,\frac{T^{\alpha }}{\alpha \Gamma (\alpha )}\right\} \right) \\
\omega &=&\mathcal{W}^{\ast }(\mathcal{W}^{\ast }e^{\mathcal{W}^{\ast }T}+1),
\end{eqnarray*}%
and $\mathcal{W}^{\ast }$ is given by (\ref{e3.20}).

\textbf{Proof.} Define the integro-differential operators $\mathcal{M}_{1}%
\mathcal{\theta }=-\mathcal{I}_{x}^{2}\mathcal{\theta }_{t}$ and $\mathcal{M}%
_{2}\phi =-\mathcal{I}_{x}^{2}\phi _{t},$ where%
\begin{equation*}
\mathcal{I}_{x}^{2}\mathcal{\theta }(x,t)=\int\limits_{0}^{x}\int%
\limits_{0}^{\xi }\mathcal{\theta }(\eta ,t)d\eta d\xi ,\text{ }\mathcal{I}%
_{x}^{2}\phi (x,t)=\int\limits_{0}^{x}\int\limits_{0}^{\xi }\phi (\eta
,t)d\eta d\xi ,
\end{equation*}%
and consider the identity%
\begin{eqnarray}
&&\left( \rho _{1}\partial _{t}^{\alpha +1}\mathcal{\theta },\mathcal{M}_{1}%
\mathcal{\theta }\right) _{L^{2}(0,L)}-\kappa _{1}((\mathcal{\theta }%
_{x}+\phi )_{x},\mathcal{M}_{1}\mathcal{\theta })_{L^{2}(0,L)}+\left(
\mathcal{\theta }_{t},\mathcal{M}_{1}\mathcal{\theta }\right) _{L^{2}(0,L)}
\notag \\
&&+(\rho _{2}\partial _{t}^{\alpha +1}\phi ,\mathcal{M}_{2}\phi
)_{L^{2}(0,L)}-\kappa _{2}((\phi _{xx},\mathcal{M}_{2}\phi
)_{L^{2}(0,L)}+\kappa _{1}((\mathcal{\theta }_{x}+\phi ),\mathcal{M}_{2}\phi
)_{L^{2}(0,L)}  \notag \\
\ &+&\left( \int\limits_{0}^{t}m(t-s)\phi _{xx}(x,s)ds,\mathcal{M}_{2}\phi
\right) _{L^{2}(0,L)}  \notag \\
&=&(F(x,t),\mathcal{M}_{1}\mathcal{\theta })_{L^{2}(0,L)}+(G(x,t),\mathcal{M}%
_{2}\phi )_{L^{2}(0,L)}.  \label{e3.2}
\end{eqnarray}%
The standard integration by parts of each term in (\ref{e3.2}) \ and
conditions (\ref{e2.3}), (\ref{e2.4}) give%
\begin{eqnarray}
\left( \rho _{1}\partial _{t}^{\alpha +1}\mathcal{\theta },\mathcal{M}_{1}%
\mathcal{\theta }\right) _{L^{2}(0,L)} &=&\frac{\rho _{1}}{2}(\partial
_{t}^{\alpha }(\mathcal{I}_{x}\mathcal{\theta }_{t}),\mathcal{I}_{x}\mathcal{%
\theta }_{t})_{L^{2}(0,L)}  \notag \\
&\geq &\frac{\rho _{1}}{2}\partial _{t}^{\alpha }\left\Vert \mathcal{I}_{x}%
\mathcal{\theta }_{t}\right\Vert _{L^{2}(0,L)},  \label{3.3*}
\end{eqnarray}%
\begin{eqnarray}
\left( \rho _{2}\partial _{t}^{\alpha +1}\phi ,\mathcal{M}_{2}\phi \right)
_{L^{2}(0,L)} &=&\frac{\rho _{2}}{2}(\partial _{t}^{\alpha }(\mathcal{I}%
_{x}\phi _{t}),\mathcal{I}_{x}\phi _{t})_{L^{2}(0,L)}  \notag \\
&\geq &\frac{\rho _{2}}{2}\partial _{t}^{\alpha }\left\Vert \mathcal{I}%
_{x}\phi _{t}\right\Vert _{L^{2}(0,L)},  \label{e3.4}
\end{eqnarray}%
\begin{equation}
-\left( \mathcal{\theta }_{t},\mathcal{M}_{1}\mathcal{\theta }\right)
_{L^{2}(0,L)}=\Vert \mathcal{I}_{x}\mathcal{\theta }_{t}\Vert
_{L^{2}(0,L)}^{2},  \label{e3.4*}
\end{equation}%
\begin{eqnarray}
\kappa _{1}(\mathcal{\theta }_{xx},\mathcal{I}_{x}^{2}\mathcal{\theta }%
_{t})_{L^{2}(0,L)} &=&\kappa _{1}\mathcal{I}_{x}^{2}\mathcal{\theta }_{t}.%
\mathcal{\theta }_{x}]_{0}^{L}-\kappa _{1}\int\limits_{0}^{L}\mathcal{I}_{x}%
\mathcal{\theta }_{t}.\mathcal{\theta }_{x}dx=\kappa _{1}\int\limits_{0}^{L}%
\mathcal{\theta }_{t}\mathcal{\theta }dx  \notag \\
&=&\frac{\kappa _{1}}{2}\frac{\partial }{\partial t}\Vert \mathcal{\theta }%
\Vert _{L^{2}(0,L)}^{2},  \label{e3.5}
\end{eqnarray}%
and in the same manner, we have%
\begin{equation}
\kappa _{2}(\phi _{xx},\mathcal{I}_{x}^{2}\phi _{t})_{L^{2}(0,L)}=\frac{%
\kappa _{2}}{2}\frac{\partial }{\partial t}\Vert \phi \Vert
_{L^{2}(0,L)}^{2},  \label{e3.6}
\end{equation}%
\begin{eqnarray}
\kappa _{1}(\phi _{x},\mathcal{I}_{x}^{2}\mathcal{\theta }%
_{t})_{L^{2}((0,L))} &=&\kappa _{1}\mathcal{I}_{x}^{2}\mathcal{\theta }%
_{t}.\phi ]_{0}^{L}dt-\kappa _{1}\int\limits_{0}^{L}\mathcal{I}_{x}\mathcal{%
\theta }_{t}.\phi dx  \notag \\
&=&-\kappa _{1}\int\limits_{0}^{L}\mathcal{I}_{x}\mathcal{\theta }_{t}.\phi
dx,  \label{e3.7}
\end{eqnarray}%
\begin{equation}
-\kappa _{1}(\mathcal{\theta }_{x},\mathcal{I}_{x}^{2}\phi
_{t})_{L^{2}(Q^{\tau })}=\kappa _{1}\int\limits_{0}^{L}\mathcal{I}_{x}\phi
_{t}.\mathcal{\theta }dx,  \label{e3.8}
\end{equation}%
\begin{eqnarray}
&&-\kappa _{1}(\phi ,\mathcal{I}_{x}^{2}\phi _{t})_{L^{2}(0,L)}  \notag \\
&=&-\kappa _{1}\mathcal{I}_{x}^{2}\phi _{t}.\Im _{x}\phi ]_{0}^{L}+\kappa
_{1}\int\limits_{0}^{\tau }\int\limits_{0}^{L}\mathcal{I}_{x}\phi _{t}.%
\mathcal{I}_{x}\phi dxdt  \notag \\
&=&\frac{1}{2}\frac{\partial }{\partial t}\Vert \mathcal{I}_{x}\phi \Vert
_{L^{2}(0,L)}^{2},  \label{e3.9}
\end{eqnarray}%
\begin{eqnarray}
&&-\left( \int\limits_{0}^{t}m(t-s).\phi _{xx}(x,s)ds,\mathcal{I}%
_{x}^{2}\phi _{t}\right) _{L^{2}(0,L)}  \notag \\
&=&-\int\limits_{0}^{L}\left( \int\limits_{0}^{t}m(t-s).\phi
_{xx}(x,s)ds\right) \mathcal{I}_{x}^{2}\phi _{t}dx  \notag \\
&=&-\int\limits_{0}^{t}m(t-s).\phi _{x}(x,s)ds).\mathcal{I}_{x}^{2}\phi
_{t}]_{0}^{L}dx+\int\limits_{0}^{L}\left( \int\limits_{0}^{t}m(t-s).\phi
_{x}(x,s)ds\right) \mathcal{I}_{x}\phi _{t}dx  \notag \\
&=&\left( \int\limits_{0}^{t}m(t-s).\phi (x,s)ds)\right) \mathcal{I}_{x}\phi
_{t}]_{0}^{L}dx-\int\limits_{0}^{L}\left( \int\limits_{0}^{t}m(t-s).\phi
(x,s)ds\right) \phi _{t}dx  \notag \\
&=&-\int\limits_{0}^{L}\left( \int\limits_{0}^{t}m(t-s).\phi (x,s)ds\right)
\phi _{t}dx.  \label{e3.10*}
\end{eqnarray}%
Substituting equalities (\ref{3.3*})-(\ref{e3.10*}) into ((\ref{e3.2})),we
obtain%
\begin{eqnarray}
&&\frac{\rho _{1}}{2}\partial _{t}^{\alpha }\left\Vert \mathcal{I}_{x}%
\mathcal{\theta }_{t}\right\Vert _{L^{2}(0,L)}+\frac{\rho _{2}}{2}\partial
_{t}^{\alpha }\left\Vert \mathcal{I}_{x}\phi _{t}\right\Vert _{L^{2}(0,L)}+%
\frac{\kappa _{1}}{2}\frac{\partial }{\partial t}\Vert \mathcal{\theta }%
\Vert _{L^{2}(0,L)}^{2}  \notag \\
&&+\frac{\kappa _{2}}{2}\frac{\partial }{\partial t}\Vert \phi \Vert
_{L^{2}(0,L)}^{2}+\Vert \mathcal{I}_{x}\mathcal{\theta }_{t}\Vert
_{L^{2}(0,L)}^{2}+\frac{1}{2}\frac{\partial }{\partial t}\Vert \mathcal{I}%
_{x}\phi \Vert _{L^{2}(0,L)}^{2}  \notag \\
&=&\kappa _{1}\int\limits_{0}^{L}\phi \mathcal{I}_{x}\mathcal{\theta }%
_{t}dx-\kappa _{1}\int\limits_{0}^{L}\mathcal{\theta I}_{x}\phi
_{t}dx-\int\limits_{0}^{L}\left( \int\limits_{0}^{t}m(t-s).\phi
(x,s)ds\right) \phi _{t}dx  \notag \\
&&-\int\limits_{0}^{L}F\mathcal{I}_{x}^{2}\mathcal{\theta }%
_{t}dx-\int\limits_{0}^{L}G\mathcal{I}_{x}^{2}\phi _{t}dx.  \label{e3.11}
\end{eqnarray}%
Replacing $t$ by $\tau $, integrating with respect to $\tau $ from zero to $%
t $ and using the given conditions, we obtain%
\begin{eqnarray}
&&\frac{\rho _{1}}{2}D^{\alpha -1}\left\Vert \mathcal{I}_{x}\mathcal{\theta }%
_{t}\right\Vert _{L^{2}(0,L)}+\frac{\rho _{2}}{2}D^{\alpha -1}\left\Vert
\mathcal{I}_{x}\phi _{t}\right\Vert _{L^{2}(0,L)}+\frac{\kappa _{1}}{2}\Vert
\mathcal{\theta }(.,t)\Vert _{L^{2}(0,L)}^{2}  \notag \\
&&+\frac{\kappa _{2}}{2}\Vert \phi (.,t)\Vert _{L^{2}(0,L)}^{2}+\Vert
\mathcal{I}_{x}\mathcal{\theta }_{\tau }\Vert _{L^{2}(0,t;L^{2}(0,L))}^{2}+%
\frac{1}{2}\Vert \mathcal{I}_{x}\phi (.,t)\Vert _{L^{2}(0,L)}^{2}  \notag \\
&=&\kappa _{1}\left( \phi ,\mathcal{I}_{x}\mathcal{\theta }_{\tau }\right)
_{L^{2}(0,t;L^{2}(0,L))}-\kappa _{1}\left( \mathcal{\theta },\mathcal{I}%
_{x}\phi _{\tau }\right) _{L^{2}(0,t;L^{2}(0,L))}  \notag \\
&&-\left( F,\mathcal{I}_{x}^{2}\mathcal{\theta }_{\tau }\right)
_{L^{2}(0,t;L^{2}(0,L))}-\left( G,\mathcal{I}_{x}^{2}\phi _{\tau }\right)
_{L^{2}(0,t;L^{2}(0,L))}  \notag \\
&&+\frac{\rho _{1}t^{1-\alpha }}{2(\Gamma (1-\alpha )(1-\alpha ))}\left\Vert
\mathcal{I}_{x}\psi \right\Vert _{L^{2}(0,L)}^{2}+\frac{\kappa _{1}}{2}\Vert
\varphi \Vert _{L^{2}(0,L)}^{2}  \notag \\
&&+\frac{\rho _{2}t^{1-\alpha }}{2(\Gamma (1-\alpha )(1-\alpha ))}\left\Vert
\mathcal{I}_{x}g\right\Vert _{L^{2}(0,L)}^{2}+\frac{\kappa _{2}}{2}\Vert
f\Vert _{L^{2}(0,L)}^{2}  \notag \\
&&+\frac{1}{2}\Vert \mathcal{I}_{x}f\Vert
_{L^{2}(0,L)}^{2}-\int\limits_{0}^{t}\int\limits_{0}^{L}\left(
\int\limits_{0}^{\tau }m(\tau -s).\phi (x,s)ds\right) \phi _{\tau }dxd\tau .
\label{e3.12*}
\end{eqnarray}%
The last term on the right-hand side needs to be evaluated as follows%
\begin{eqnarray}
&&-\int\limits_{0}^{t}\int\limits_{0}^{L}\left( \int\limits_{0}^{\tau
}m(\tau -s).\phi (x,s)ds\right) \phi _{\tau }dxd\tau  \notag \\
&=&-\int\limits_{0}^{L}\left( \int\limits_{0}^{\tau }m(\tau -s).\phi
(x,s)ds\right) \phi dx]_{0}^{t}+\int\limits_{0}^{\tau
}\int\limits_{0}^{L}m(0)\phi ^{2}dxd\tau  \notag \\
&&+\int\limits_{0}^{t}\int\limits_{0}^{L}\left( \int\limits_{0}^{\tau
}m^{\prime }(\tau -s).\phi (x,s)ds\right) \phi (x,\tau )dxd\tau  \notag \\
&=&-\int\limits_{0}^{L}\left( \int\limits_{0}^{t}m(t-s).\phi (x,s)ds\right)
\phi (x,t)dx+m(0)\Vert \phi \Vert _{L^{2}(0,t;L^{2}(0,L))}^{2}  \notag \\
&&+\int\limits_{0}^{t}\int\limits_{0}^{L}\left( \int\limits_{0}^{\tau
}m^{\prime }(\tau -s).\phi (x,s)ds\right) \phi dxd\tau .  \label{e3.13*}
\end{eqnarray}%
By replacing (\ref{e3.13*}) into (\ref{e3.12*}), and estimating different
terms on the right-hand side of (by using Cauchy $\epsilon $ inequality, a
Poincare type inequality) (\ref{e3.12*}) as follows%
\begin{eqnarray}
&&\kappa _{1}\left( \phi ,\mathcal{I}_{x}\mathcal{\theta }_{t}\right)
_{L^{2}(0,t;L^{2}(0,L))}  \notag \\
&\leq &\frac{\kappa _{1}\epsilon _{1}}{2}\Vert \phi \Vert
_{L^{2}(0,t;L^{2}(0,L))}^{2}+\frac{\kappa _{1}}{2\epsilon _{1}}\Vert
\mathcal{I}_{x}\mathcal{\theta }_{\tau }\Vert _{L^{2}(0,t;L^{2}(0,L))}^{2},\
\ \ \ \ \ \ \ \ \ \   \label{e3.12}
\end{eqnarray}%
\begin{eqnarray}
&&-\kappa _{1}\left( \mathcal{\theta },\mathcal{I}_{x}\phi _{\tau }\right)
_{L^{2}(0,t;L^{2}(0,L))}  \notag \\
&\leq &\frac{\kappa _{1}\epsilon _{2}}{2}\Vert \mathcal{\theta }\Vert
_{L^{2}(0,t;L^{2}(0,L))}^{2}+\frac{\kappa _{1}}{2\epsilon _{2}}\Vert
\mathcal{I}_{x}\phi _{\tau }\Vert _{L^{2}(0,t;L^{2}(0,L))}^{2},\ \ \ \ \ \ \
\ \   \label{e3.13}
\end{eqnarray}%
\begin{eqnarray}
&&-\int\limits_{0}^{L}\left( \int\limits_{0}^{t}m(t-s).\phi (x,s)ds\right)
\phi (x,t)dx  \notag \\
&\leq &\frac{\epsilon _{3}}{2}\int\limits_{0}^{L}\phi ^{2}(x,t)dx+\frac{1}{%
2\epsilon _{3}}\int\limits_{0}^{L}\left( \int\limits_{0}^{t}m(t-s)\phi
(x,s)ds\right) ^{2}dx  \notag \\
&\leq &\frac{\epsilon _{3}}{2}\int\limits_{0}^{L}\phi ^{2}(x,t)dx+\frac{1}{%
2\epsilon _{3}}\int\limits_{0}^{L}\left(
\int\limits_{0}^{t}m^{2}(t-s)ds\right) \left( \int\limits_{0}^{t}\phi
^{2}(x,s)ds\right) dx  \notag \\
&\leq &\frac{\epsilon _{3}}{2}\int\limits_{0}^{L}\phi ^{2}(x,t)dx+\frac{T}{%
2\epsilon _{3}}\sup\limits_{0\leq t\leq
T}m^{2}(t)\int\limits_{0}^{L}\int\limits_{0}^{t}\phi ^{2}dxd\tau ,
\label{e3.14}
\end{eqnarray}%
\begin{eqnarray}
&&\int\limits_{0}^{t}\int\limits_{0}^{L}\left( \int\limits_{0}^{\tau
}m^{\prime }(\tau -s).\phi (x,s)ds\right) \phi dxd\tau  \notag \\
&\leq &\frac{\epsilon _{4}}{2}\int\limits_{0}^{t}\int\limits_{0}^{L}\phi
^{2}dxd\tau +\frac{1}{2\epsilon _{4}}\int\limits_{0}^{t}\int\limits_{0}^{L}%
\left( \int\limits_{0}^{\tau }m^{\prime 2}(\tau -s)ds\right) \left(
\int\limits_{0}^{\tau }\phi ^{2}(x,s)ds\right) dxd\tau  \notag \\
&\leq &\frac{\epsilon _{4}}{2}\int\limits_{0}^{t}\int\limits_{0}^{L}\phi
^{2}dxd\tau +\frac{T}{2\epsilon _{4}}\sup\limits_{0\leq t\leq T}m^{\prime
2}(t)\int\limits_{0}^{t}\int\limits_{0}^{L}\left( \int\limits_{0}^{\tau
}\phi ^{2}(x,s)ds\right) dxd\tau  \notag \\
&=&\frac{\epsilon _{4}}{2}\int\limits_{0}^{t}\int\limits_{0}^{L}\phi
^{2}dxd\tau +\frac{T}{2\epsilon _{4}}\sup\limits_{0\leq t\leq T}m^{\prime
2}(t)\int\limits_{0}^{L}\left[ \left( \tau \int\limits_{0}^{\tau }\phi
^{2}(x,s)ds\right) _{0}^{t}-\int\limits_{0}^{t}\tau \phi ^{2}d\tau \right] dx
\notag \\
&=&\frac{\epsilon _{4}}{2}\int\limits_{0}^{t}\int\limits_{0}^{L}\phi
^{2}dxd\tau +\frac{T}{2\epsilon _{4}}\sup\limits_{0\leq t\leq T}m^{\prime
2}(t)\int\limits_{0}^{L}\left[ \int\limits_{0}^{t}(t-\tau )\phi ^{2}(x,\tau
)d\tau \right] dx  \notag \\
&\leq &\frac{\epsilon _{4}}{2}\int\limits_{0}^{t}\int\limits_{0}^{L}\phi
^{2}dxd\tau +\frac{T^{2}}{2\epsilon _{4}}\sup\limits_{0\leq t\leq
T}m^{\prime 2}(t)\int\limits_{0}^{L}\int\limits_{0}^{t}\phi ^{2}d\tau dx
\notag \\
&=&\left( \frac{\epsilon _{4}}{2}+\frac{T^{2}}{2\epsilon _{4}}%
\sup\limits_{0\leq t\leq T}m^{\prime 2}(t)\right)
\int\limits_{0}^{t}\int\limits_{0}^{L}\phi ^{2}dxd\tau ,\ \   \label{e3.15}
\end{eqnarray}%
\begin{equation}
-\left( F,\mathcal{I}_{x}^{2}\mathcal{\theta }_{\tau }\right)
_{L^{2}(0,t;L^{2}(0,L))}\leq \frac{\epsilon _{5}}{2}\Vert F\Vert
_{L^{2}(0,t;L^{2}(0,L))}^{2}+\frac{L^{2}}{2\epsilon _{5}}\Vert \mathcal{I}%
_{x}\mathcal{\theta }_{\tau }\Vert _{L^{2}(0,t;L^{2}(0,L))}^{2},
\label{e3.16*}
\end{equation}%
\begin{equation}
-\left( G,\mathcal{I}_{x}^{2}\phi _{\tau }\right)
_{L^{2}(0,t;L^{2}(0,L))}\leq \frac{\epsilon _{6}}{2}\Vert G\Vert
_{L^{2}(0,t;L^{2}(0,L))}^{2}+\frac{L^{2}}{4\epsilon _{6}}\Vert \mathcal{I}%
_{x}\phi _{\tau }\Vert _{L^{2}(0,t;L^{2}(0,L)),}^{2}  \label{e3.17*}
\end{equation}%
Combination of (\ref{e3.12})-(\ref{e3.17*}) and (\ref{e3.12*}), leads to%
\begin{eqnarray}
&&\frac{\rho _{1}}{2}D^{\alpha -1}\left\Vert \mathcal{I}_{x}\mathcal{\theta }%
_{t}\right\Vert _{L^{2}(0,L)}+\frac{\rho _{2}}{2}D^{\alpha -1}\left\Vert
\mathcal{I}_{x}\phi _{t}\right\Vert _{L^{2}(0,L)}+\frac{\kappa _{1}}{2}\Vert
\mathcal{\theta }(.,t)\Vert _{L^{2}(0,L)}^{2}  \notag \\
&&+\frac{\kappa _{2}}{2}\Vert \phi (.,t)\Vert _{L^{2}(0,L)}^{2}+\Vert
\mathcal{I}_{x}\mathcal{\theta }_{\tau }\Vert _{L^{2}(0,t;L^{2}(0,L))}^{2}+%
\frac{1}{2}\Vert \mathcal{I}_{x}\phi (.,t)\Vert _{L^{2}(0,L)}^{2}  \notag \\
&\leq &\frac{\rho _{1}T^{1-\alpha }L^{2}}{4(\Gamma (1-\alpha )(1-\alpha ))}%
\left\Vert \psi \right\Vert _{L^{2}(0,L)}^{2}+\frac{\kappa _{1}}{2}\Vert
\varphi \Vert _{L^{2}(0,L)}^{2}  \notag \\
&&+\frac{\rho _{2}T^{1-\alpha }L^{2}}{4(\Gamma (1-\alpha )(1-\alpha ))}%
\left\Vert g\right\Vert _{L^{2}(0,L)}^{2}+\left( \frac{\kappa _{2}}{2}+\frac{%
L^{2}}{4}\right) \Vert f\Vert _{L^{2}(0,L)}^{2}  \notag \\
&&+\left( \frac{\kappa _{1}\epsilon _{1}}{2}+\frac{T}{2\epsilon _{3}}%
\sup\limits_{0\leq t\leq T}m^{2}(t)+\frac{\epsilon _{4}}{2}+\frac{T^{2}}{%
2\epsilon _{4}}\sup\limits_{0\leq t\leq T}m^{\prime 2}(t)+m(0)\right) \Vert
\phi \Vert _{L^{2}(0,t;L^{2}(0,L))}^{2}  \notag \\
&&+\left( \frac{\kappa _{1}}{2\epsilon _{1}}+\frac{L^{2}}{2\epsilon _{5}}%
\right) \Vert \mathcal{I}_{x}\mathcal{\theta }_{\tau }\Vert
_{L^{2}(0,t;L^{2}(0,L))}^{2}+\left( \frac{\kappa _{1}}{2\epsilon _{2}}+\frac{%
L^{2}}{4\epsilon _{6}}\right) \Vert \mathcal{I}_{x}\phi _{\tau }\Vert
_{L^{2}(0,t;L^{2}(0,L))}^{2}  \notag \\
&&+\frac{\kappa _{1}\epsilon _{2}}{2}\Vert \mathcal{\theta }\Vert
_{L^{2}(0,t;L^{2}(0,L))}^{2}+\frac{\epsilon _{3}}{2}\Vert \phi (.,t)\Vert
_{L^{2}(0,L)}^{2}+\frac{\epsilon _{5}}{2}\Vert F\Vert
_{L^{2}(0,t;L^{2}(0,L))}^{2}  \notag \\
&&+\frac{\epsilon _{6}}{2}\Vert G\Vert _{L^{2}(0,t;L^{2}(0,L))}^{2}.
\label{e3.18}
\end{eqnarray}%
The choice $\varepsilon _{1}=\kappa _{1},$ $\varepsilon _{5}=L^{2}/2,$ $%
\varepsilon _{3}=\kappa _{2}/2,$ $\varepsilon _{2}=\varepsilon
_{4}=\varepsilon _{6}=1,$ and cancellation of the last term on the left-hand
side of (\ref{e3.18}) reduce it to the following estimate%
\begin{eqnarray}
&&D^{\alpha -1}\left\Vert \mathcal{I}_{x}\mathcal{\theta }_{t}\right\Vert
_{L^{2}(0,L)}+D^{\alpha -1}\left\Vert \mathcal{I}_{x}\phi _{t}\right\Vert
_{L^{2}(0,L)}+\Vert \mathcal{\theta }(.,t)\Vert _{L^{2}(0,L)}^{2}+\Vert \phi
(.,t)\Vert _{L^{2}(0,L)}^{2}  \notag \\
&\leq &\mathcal{W}^{\ast }\left( \Vert \mathcal{I}_{x}\mathcal{\theta }%
_{\tau }\Vert _{L^{2}(0,t;L^{2}(0,L))}^{2}+\Vert \mathcal{I}_{x}\phi _{\tau
}\Vert _{L^{2}(0,t;L^{2}(0,L))}^{2}+\Vert \mathcal{\theta }\Vert
_{L^{2}(0,t;L^{2}(0,L))}^{2}+\Vert \phi \Vert
_{L^{2}(0,t;L^{2}(0,L))}^{2}\right.  \notag \\
&&+\Vert \varphi \Vert _{L^{2}(0,L)}^{2}+\left\Vert \psi \right\Vert
_{L^{2}(0,L)}+\Vert f\Vert _{L^{2}(0,L)}^{2}+\Vert g\Vert
_{L^{2}(0,L)}^{2}+\Vert F\Vert _{L^{2}(0,t;L^{2}(0,L))}^{2}  \label{e3.19*}
\\
&&\left. +\Vert G\Vert _{L^{2}(0,t;L^{2}(0,L))}^{2}\right) ,  \notag
\end{eqnarray}%
where%
\begin{eqnarray}
\mathcal{W}^{\ast } &=&\max \left( \frac{\kappa _{1}^{2}}{2}+\frac{T}{\kappa
_{2}}\sup\limits_{0\leq t\leq T}m^{2}(t)+\frac{1}{2}+\frac{T^{2}}{2}%
\sup\limits_{0\leq t\leq T}m^{\prime 2}(t)+m(0),\text{ }\frac{3}{2},\text{ }%
\frac{\kappa _{1}}{2}+\frac{L^{2}}{4},\frac{\kappa _{2}}{4},\right.  \notag
\\
&&\left. \frac{\rho _{1}T^{1-\alpha }L^{2}}{4(\Gamma (1-\alpha )(1-\alpha ))}%
,\frac{\rho _{2}T^{1-\alpha }L^{2}}{4(\Gamma (1-\alpha )(1-\alpha ))}\right)
/\min \left( \frac{\rho _{1}}{2},\frac{\rho _{2}}{2},\frac{\kappa _{1}}{2},%
\frac{\kappa _{2}}{2},\frac{1}{2}\right) .  \label{e3.20}
\end{eqnarray}%
By omitting the first and second term on the left-hand side of (\ref{e3.19*}%
), and applying the Gronwall-Bellman lemma ( [46]), where%
\begin{equation}
\left\{
\begin{array}{c}
E(t)=\Vert \mathcal{\theta }\Vert _{L^{2}(0,t;L^{2}(0,L))}^{2}+\Vert \phi
\Vert _{L^{2}(0,t;L^{2}(0,L))}^{2} \\
\frac{dE}{dt}=\Vert \mathcal{\theta }(.,t)\Vert _{L^{2}(0,L)}^{2}+\Vert \phi
(.,t)\Vert _{L^{2}(0,L)}^{2}, \\
\mathcal{Q}(0)=0.%
\end{array}%
\right.  \label{e3.21*}
\end{equation}%
we obtain%
\begin{eqnarray}
y(t) &\leq &\mathcal{W}^{\ast }e^{\mathcal{W}^{\ast }t}\left( \Vert \mathcal{%
I}_{x}\mathcal{\theta }_{\tau }\Vert _{L^{2}(0,t;L^{2}(0,L))}^{2}+\Vert
\mathcal{I}_{x}\phi _{\tau }\Vert _{L^{2}(0,t;L^{2}(0,L))}^{2}\right.  \notag
\\
&&+\Vert \varphi \Vert _{L^{2}(0,L)}^{2}+\left\Vert \psi \right\Vert
_{L^{2}(0,L)}+\Vert f\Vert _{L^{2}(0,L)}^{2}+\Vert g\Vert _{L^{2}(0,L)}^{2}
\label{e3.22**} \\
&&+\left. \Vert F\Vert _{L^{2}(0,t;L^{2}(0,L))}^{2}+\Vert G\Vert
_{L^{2}(0,t;L^{2}(0,L))}^{2}\right) .  \notag
\end{eqnarray}%
Then by omitting the last two terms on the left-hand side of (\ref{e3.19*}),
and using (\ref{e3.22**}), we have%
\begin{eqnarray}
&&D^{\alpha -1}\left( \left\Vert \mathcal{I}_{x}\mathcal{\theta }%
_{t}\right\Vert _{L^{2}(0,L)}+\left\Vert \mathcal{I}_{x}\phi _{t}\right\Vert
_{L^{2}(0,L)}\right)  \notag \\
&\leq &\mathcal{W}^{\ast }(\mathcal{W}^{\ast }e^{\mathcal{W}^{\ast
}T}+1)\left( \Vert \mathcal{I}_{x}\mathcal{\theta }_{\tau }\Vert
_{L^{2}(0,t;L^{2}(0,L))}^{2}+\Vert \mathcal{I}_{x}\phi _{\tau }\Vert
_{L^{2}(0,t;L^{2}(0,L))}^{2}\right.  \notag \\
&&+\Vert \varphi \Vert _{L^{2}(0,L)}^{2}+\left\Vert \psi \right\Vert
_{L^{2}(0,L)}+\Vert f\Vert _{L^{2}(0,L)}^{2}+\Vert g\Vert _{L^{2}(0,L)}^{2}
\label{e3.23*} \\
&&+\left. \Vert F\Vert _{L^{2}(0,t;L^{2}(0,L))}^{2}+\Vert G\Vert
_{L^{2}(0,t;L^{2}(0,L))}^{2}\right) .  \notag
\end{eqnarray}%
Now, Lemma \ 3.2, can be applied to remove the first two terms on the
right-hand side of (\ref{e3.23*}), by taking%
\begin{equation}
\left\{
\begin{array}{c}
\mathcal{Q}(t)=\Vert \mathcal{I}_{x}\mathcal{\theta }_{\tau }\Vert
_{L^{2}(0,t;L^{2}(0,L))}^{2}+\Vert \mathcal{I}_{x}\phi _{\tau }\Vert
_{L^{2}(0,t;L^{2}(0,L))}^{2} \\
^{C}\partial _{t}^{\beta }\mathcal{Q}(t)=D^{\alpha -1}\left( \left\Vert
\mathcal{I}_{x}\mathcal{\theta }_{t}\right\Vert _{L^{2}(0,L)}+\left\Vert
\mathcal{I}_{x}\phi _{t}\right\Vert _{L^{2}(0,L)}\right) \\
\mathcal{Q}(0)=0.%
\end{array}%
\right.  \label{e3.24*}
\end{equation}%
From (\ref{e3.23*}), it follows that%
\begin{eqnarray}
\mathcal{Q}(t) &\leq &\mathcal{M}\left\{ D^{-1-\alpha }\left( \Vert F\Vert
_{L^{2}(0,L)}^{2}+\Vert G\Vert _{L^{2}(0,L)}^{2}\right) \right.  \notag \\
&&+\left. \Vert \varphi \Vert _{L^{2}(0,L)}^{2}+\left\Vert \psi \right\Vert
_{L^{2}(0,L)}+\Vert f\Vert _{L^{2}(0,L)}^{2}+\Vert g\Vert
_{L^{2}(0,L)}^{2}\right\} ,  \label{e3.25*}
\end{eqnarray}%
where%
\begin{equation*}
\mathcal{M}=\Gamma (\alpha )E_{\alpha ,\alpha }(\omega t^{\alpha })\left(
\max \left\{ 1,\frac{T^{\alpha }}{\alpha \Gamma (\alpha )}\right\} \right) ,
\end{equation*}%
with
\begin{equation*}
\omega =\mathcal{W}^{\ast }(\mathcal{W}^{\ast }e^{\mathcal{W}^{\ast }T}+1).
\end{equation*}%
Now since%
\begin{equation}
D^{-1-\alpha }\left( \Vert F\Vert _{L^{2}(0,L)}^{2}+\Vert G\Vert
_{L^{2}(0,L)}^{2}\right) \leq \frac{T^{\alpha }}{\Gamma (\alpha +1)}%
\int\limits_{0}^{t}\left( \Vert F\Vert _{L^{2}(0,L)}^{2}+\Vert G\Vert
_{L^{2}(0,L)}^{2}\right) d\tau ,  \label{e3.26*}
\end{equation}%
then, we infer from (\ref{e3.25*}), (\ref{e3.26*}), and (\ref{e3.19*}) the
following inequality%
\begin{eqnarray}
&&D^{\alpha -1}\left\Vert \mathcal{I}_{x}\mathcal{\theta }_{t}\right\Vert
_{L^{2}(0,L)}+D^{\alpha -1}\left\Vert \mathcal{I}_{x}\phi _{t}\right\Vert
_{L^{2}(0,L)}+\Vert \mathcal{\theta }(.,t)\Vert _{L^{2}(0,L)}^{2}+\Vert \phi
(.,t)\Vert _{L^{2}(0,L)}^{2}  \notag \\
&\leq &\mathcal{F}^{\ast }\left( \Vert \psi \Vert _{L^{2}(0,L)}^{2}+\Vert
\varphi \Vert _{L^{2}(0,L)}^{2}+\Vert f\Vert _{L^{2}(0,L)}^{2}+\Vert g\Vert
_{L^{2}(0,L)}^{2}\right.  \notag \\
&&\left. +\Vert F\Vert _{L^{2}(0,T;L^{2}(0,L))}^{2}+\Vert G\Vert
_{L^{2}(0,T;L^{2}(0,L))}^{2}\right) .  \label{e3.27*}
\end{eqnarray}%
The first estimate (\ref{e3.1}) follows, if we disregard the first and
second term on the left-hand side of (\ref{e3.27*}), and pass to the
supremum with respect to $t$ over $(0,T).$ The second estimate (\ref{e3.1**}%
) follows from (\ref{e3.27*}) if we drop the last two terms on the left-hand
side of the inequality (\ref{e3.27*}).

Since the range of the operator $\mathcal{G}$\ is subset of $\mathcal{E}$,
that is $R\left( \mathcal{G}\right) \subset \mathcal{E}$, so we extend $%
\mathcal{G}$ so that inequality (\ref{e3.27*}) holds for the extension, and $%
R\left( \overline{\mathcal{G}}\right) =\mathcal{E}$. We can easily show that
the following

\textbf{Proposition 3.2}. The unbounded operator $\mathcal{G}:\mathcal{B}%
\rightarrow \mathcal{E}$ admits a closure $\overline{\mathcal{G}\text{ }}$
with domain of definition $D(\overline{\mathcal{G}\text{ }}).$

\textbf{Definition 3.3. }The solution of the equation $\overline{\mathcal{G}%
\text{ }}\mathcal{Z=}H=(\{F,\varphi ,\psi \},\{G,f,g\})$ is called a strong
solution of problem (\ref{e2.1}), (\ref{e2.3}), (\ref{e2.4}).

The a priori estimate (\ref{e3.1}), can be extended to%
\begin{equation}
\Vert \mathcal{Z}\Vert _{\mathcal{B}}^{2}\leq \mathcal{F}^{\ast }\Vert
\overline{\mathcal{G}\text{ }}\mathcal{Z}\Vert _{\mathcal{E}}^{2},\text{ \ \
}\forall \mathcal{Z\in }D(\overline{\mathcal{G}\text{ }}).  \label{e3.28*}
\end{equation}%
The estimate (\ref{e3.28*}) shows that the operator $\overline{\mathcal{G}%
\text{ }}$ is one to one and that $\overline{\mathcal{G}\text{ }}^{-1}$ is
continuous from $R(\overline{\mathcal{G}\text{ }})$ onto $\mathcal{B}.$
Consequently if a strong solution of problem (\ref{e2.1}), (\ref{e2.3}), (%
\ref{e2.4}) exists, it is unique and depends continuously on the input data $%
\ \varphi ,\psi ,f,g$ and the external forces $F$ and $G.$ Also as a
consequence of (\ref{e3.28*}), the set $R(\overline{\mathcal{G}\text{ }}%
)\subset \mathcal{E}$ is closed and $R(\overline{\mathcal{G}\text{ }})=%
\overline{R(\mathcal{G})}.$

\section{Existence of solution}

\textbf{Theorem 4.1}. Problem (\ref{e2.1}), (\ref{e2.3}), (\ref{e2.4})
admits a unique strong solution $\mathcal{Z}=\left( \overline{\mathcal{G}%
\text{ }}\right) ^{-1}(\{F,\varphi ,\psi \},\{G,f,g\})=\overline{\mathcal{G}%
^{-1}}(\{F,\varphi ,\psi \},\{G,f,g\}),$which depend continuously on the
given data, for all $F,G\in L^{2}(0,T;L^{2}(0,L)),$and $\varphi ,\psi
,f,g\in L^{2}(0,L).$

\textbf{Proof}. It follows from above that in order to prove the existence
of the strongly generalized solution of problem (\ref{e2.1}), (\ref{e2.3}), (%
\ref{e2.4}), it suffices to prove that $\overline{R(\mathcal{G})}=\mathcal{E}
$. To this end, we first prove the density in the following special case.

\textbf{Theorem 4.2}. If for some function $W(x,t)=(\Lambda
_{1}(x,t),\Lambda _{2}(x,t))\in (L^{2}(0,T;L^{2}(0,L)))^{2}$ and for
elements $\ \mathcal{Z}\in D_{0}(\mathcal{G})=\{\mathcal{Z}\ :\ \mathcal{Z}%
\in D(\mathcal{G})\ $and$\ \Gamma _{i}\mathcal{\theta }=\Gamma _{i}\phi =0\
,\ i=1,2\ \}$ we have
\begin{equation}
(\mathcal{S}_{1}(\mathcal{\theta },\phi ),\Lambda
_{1})_{L^{2}(0,T;L^{2}(0,L))}+(\mathcal{S}_{2}(\mathcal{\theta },\phi
),\Lambda _{2})_{L^{2}(0,T;L^{2}(0,L))}=0,  \label{e4.1}
\end{equation}%
then $W(x,t)=(\Lambda _{1}(x,t),\Lambda _{2}(x,t))=(0,0)$ a.e in $Q^{T}$.

\textbf{Proof.} The identity (\ref{e4.1}) is equivalent to%
\begin{eqnarray}
&&\int\limits_{0}^{T}(\rho _{1}\partial _{t}^{\alpha +1}\mathcal{\theta }%
,\Lambda _{1})_{L^{2}(0,L)}dt-\kappa _{1}\int\limits_{0}^{T}(\mathcal{\theta
}_{xx},\Lambda _{1})_{L^{2}(0,L)}dt-\kappa _{1}\int\limits_{0}^{T}(\phi
_{x},\Lambda _{1})_{L^{2}(0,L)}dt+\int\limits_{0}^{T}(\mathcal{\theta }%
_{t},\Lambda _{1})_{L^{2}(0,L)}dt  \notag \\
\ &&+\int\limits_{0}^{T}(\rho _{2}\partial _{t}^{\alpha +1}\phi ,\Lambda
_{2})_{L^{2}(0,L)}dt-\kappa _{2}\int\limits_{0}^{T}(\phi _{xx},\Lambda
_{2})_{L^{2}(0,L)}dt+\kappa _{1}\int\limits_{0}^{T}((\mathcal{\theta }%
_{x},\Lambda _{2})_{L^{2}(0,L)}dt  \notag \\
\ &&+\kappa _{1}\int\limits_{0}^{T}(\phi ,\Lambda
_{2})_{L^{2}(0,L)}dt+\int\limits_{0}^{T}(\int\limits_{0}^{t}m(t-s)\phi
_{xx}(x,s)ds,\Lambda _{2})_{L^{2}(0,L)}dsdt  \label{e4.2} \\
&=&0.\   \notag
\end{eqnarray}%
Assume that the functions $\xi (x,t),$ $\eta (x,t)$and satisfy the boundary
and the initial conditions (\ref{e2.3}), and (\ref{e2.4}) and such that $\xi
,\eta ,$ $\xi _{x},\eta _{x},$ $\mathcal{I}_{t}\xi ,$ $\mathcal{I}_{t}\eta ,$
$\mathcal{I}_{t}\mathcal{I}_{x}^{2}\xi ,$ $\mathcal{I}_{t}\mathcal{I}%
_{x}^{2}\eta ,$ $\mathcal{I}_{t}^{2}\xi ,$ $\mathcal{I}_{t}^{2}\eta $ and $%
\partial _{t}^{\beta +1}\xi ,$ $\partial _{t}^{\beta +1}\eta \in
L^{2}(0,T;L^{2}(0,L))$, we then set%
\begin{equation}
\mathcal{\theta }(x,t)=\mathcal{I}_{t}^{2}\xi
=\int\limits_{0}^{t}\int\limits_{0}^{s}\xi (x,z)dzds,\text{ }\phi (x,t)=%
\mathcal{I}_{t}^{2}\eta =\int\limits_{0}^{t}\int\limits_{0}^{s}\eta
(x,z)dzds,  \label{e4.3}
\end{equation}%
and introduce the functions
\begin{equation}
\Lambda _{1}(x,t)=\mathcal{I}_{t}\xi +\mathcal{I}_{x}^{2}\mathcal{I}_{t}\xi ,%
\text{ }\Lambda _{2}(x,t)=\mathcal{I}_{t}\eta +\mathcal{I}_{x}^{2}\mathcal{I}%
_{t}\eta .  \label{e4.3*}
\end{equation}%
Equation (\ref{e4.2}) then reduces to%
\begin{eqnarray}
&&\int\limits_{0}^{T}(\rho _{1}\partial _{t}^{\alpha +1}\mathcal{I}%
_{t}^{2}\xi ,\mathcal{I}_{t}\xi +\mathcal{I}_{x}^{2}\mathcal{I}_{t}\xi
)_{L^{2}(0,L)}dt-\kappa _{1}\int\limits_{0}^{T}(\mathcal{I}_{t}^{2}\xi _{xx},%
\mathcal{I}_{t}\xi +\mathcal{I}_{x}^{2}\mathcal{I}_{t}\xi )_{L^{2}(0,L)}dt
\notag \\
&&-\kappa _{1}\int\limits_{0}^{T}(\mathcal{I}_{t}^{2}\eta _{x},\mathcal{I}%
_{t}\xi +\mathcal{I}_{x}^{2}\mathcal{I}_{t}\xi
)_{L^{2}(0,L)}dt+\int\limits_{0}^{T}(\mathcal{I}_{t}\xi ,\mathcal{I}_{t}\xi +%
\mathcal{I}_{x}^{2}\mathcal{I}_{t}\xi )_{L^{2}(0,L)}dt  \notag \\
&&+\int\limits_{0}^{T}(\rho _{2}\partial _{t}^{\alpha +1}\mathcal{I}%
_{t}^{2}\eta ,\mathcal{I}_{t}\eta +\mathcal{I}_{x}^{2}\mathcal{I}_{t}\eta
)_{L^{2}(0,L)}dt-\kappa _{2}\int\limits_{0}^{T}(\mathcal{I}_{t}^{2}\eta
_{xx},\mathcal{I}_{t}\eta +\mathcal{I}_{x}^{2}\mathcal{I}_{t}\eta
)_{L^{2}(0,L)}dt  \notag \\
&&+\kappa _{1}\int\limits_{0}^{T}((\mathcal{I}_{t}^{2}\xi _{x},\mathcal{I}%
_{t}\eta +\mathcal{I}_{x}^{2}\mathcal{I}_{t}\eta )_{L^{2}(0,L)}dt+\kappa
_{1}\int\limits_{0}^{T}(\mathcal{I}_{t}^{2}\eta ,\mathcal{I}_{t}\eta +%
\mathcal{I}_{x}^{2}\mathcal{I}_{t}\eta )_{L^{2}(0,L)}dt  \notag \\
&&+\int\limits_{0}^{T}(\int\limits_{0}^{t}m(t-s)\mathcal{I}_{s}^{2}\eta
_{xx}(x,s)ds,\mathcal{I}_{t}\eta +\mathcal{I}_{x}^{2}\mathcal{I}_{t}\eta
)_{L^{2}(0,L)}dt  \notag \\
&=&0.  \label{e4.4}
\end{eqnarray}%
Invoking boundary integral conditions and carrying out appropriate
integrations by parts of each inner product term, we have%
\begin{eqnarray}
&&(\rho _{1}\partial _{t}^{\alpha +1}\mathcal{I}_{t}^{2}\xi ,\mathcal{I}%
_{t}\xi +\mathcal{I}_{x}^{2}\mathcal{I}_{t}\xi )_{L^{2}(0,L)}  \notag \\
&=&(\rho _{1}\partial _{t}^{\alpha }\mathcal{I}_{t}\xi ,\mathcal{I}_{t}\xi
)+(\rho _{1}\partial _{t}^{\alpha }\mathcal{I}_{x}\mathcal{I}_{t}\xi ,%
\mathcal{I}_{x}\mathcal{I}_{t}\xi )_{L^{2}(0,L)},  \label{e4.5}
\end{eqnarray}%
\begin{eqnarray}
&&-\kappa _{1}(\mathcal{I}_{t}^{2}\xi _{xx},\mathcal{I}_{t}\xi +\mathcal{I}%
_{x}^{2}\mathcal{I}_{t}\xi )_{L^{2}(0,L)}  \notag \\
&=&\frac{\kappa _{1}\partial }{2\partial t}\Vert \mathcal{I}_{t}^{2}\xi
_{x}\Vert _{L^{2}(0,L)}^{2}+\frac{\kappa _{1}\partial }{2\partial t}\Vert
\mathcal{I}_{t}^{2}\xi \Vert _{L^{2}(0,L)}^{2},  \label{e4.6}
\end{eqnarray}%
\begin{eqnarray}
&&-\kappa _{1}(\mathcal{I}_{t}^{2}\eta _{x},\mathcal{I}_{t}\xi +\mathcal{I}%
_{x}^{2}\mathcal{I}_{t}\xi )_{L^{2}(0,L)}  \notag \\
&=&-\kappa _{1}(\mathcal{I}_{t}^{2}\eta _{x},\mathcal{I}_{t}\xi
)_{L^{2}(0,L)}+\kappa _{1}(\mathcal{I}_{t}^{2}\eta ,\mathcal{I}_{x}\mathcal{I%
}_{t}\xi )_{L^{2}(0,L)},  \label{e4.7}
\end{eqnarray}%
\begin{eqnarray}
&&(\mathcal{I}_{t}\xi ,\mathcal{I}_{t}\xi +\mathcal{I}_{x}^{2}\mathcal{I}%
_{t}\xi )_{L^{2}(0,L)}  \notag \\
&=&\Vert \Im _{t}\xi \Vert _{L^{2}(0,L)}^{2}-\Vert \mathcal{I}_{x}\Im
_{t}\xi \Vert _{L^{2}(0,L)}^{2},  \label{e4.8*}
\end{eqnarray}%
\begin{eqnarray}
&&(\rho _{2}\partial _{t}^{\alpha +1}\mathcal{I}_{t}^{2}\eta ,\mathcal{I}%
_{t}\eta +\mathcal{I}_{x}^{2}\mathcal{I}_{t}\eta )_{L^{2}(0,L)}  \notag \\
&=&(\rho _{2}\partial _{t}^{\alpha }\mathcal{I}_{t}\eta ,\mathcal{I}_{t}\eta
)+(\rho _{2}\partial _{t}^{\alpha }\mathcal{I}_{x}\mathcal{I}_{t}\eta ,%
\mathcal{I}_{x}\mathcal{I}_{t}\eta )_{L^{2}(0,L)},  \label{e4.9}
\end{eqnarray}%
\begin{eqnarray}
&&-\kappa _{2}(\mathcal{I}_{t}^{2}\eta _{xx},\mathcal{I}_{t}\eta +\mathcal{I}%
_{x}^{2}\mathcal{I}_{t}\eta )_{L^{2}(0,L)}  \notag \\
&=&\frac{\kappa _{2}\partial }{2\partial t}\Vert \mathcal{I}_{t}^{2}\eta
_{x}\Vert _{L^{2}(0,L)}^{2}+\frac{\kappa _{2}\partial }{2\partial t}\Vert
\mathcal{I}_{t}^{2}\eta \Vert _{L^{2}(0,L)}^{2}  \label{e4.10*}
\end{eqnarray}%
\begin{eqnarray}
&&\kappa _{1}(\mathcal{I}_{t}^{2}\xi _{x},\mathcal{I}_{t}\eta +\mathcal{I}%
_{x}^{2}\mathcal{I}_{t}\eta )_{L^{2}(0,L)}  \notag \\
&=&\kappa _{1}(\mathcal{I}_{t}^{2}\xi _{x},\mathcal{I}_{t}\eta
)_{L^{2}(0,L)}-\kappa _{1}(\mathcal{I}_{t}^{2}\xi ,\mathcal{I}_{x}\mathcal{I}%
_{t}\eta )_{L^{2}(0,L)},  \label{e4.11}
\end{eqnarray}%
\begin{eqnarray}
&&(\int\limits_{0}^{t}m(t-s)\mathcal{I}_{s}^{2}\eta _{xx}(x,s)ds,\mathcal{I}%
_{t}\eta +\mathcal{I}_{x}^{2}\mathcal{I}_{t}\eta )_{L^{2}(0,L)}  \notag \\
&=&-(\int\limits_{0}^{t}m(t-s)\mathcal{I}_{s}^{2}\eta _{x}(x,s)ds,\mathcal{I}%
_{t}\eta _{x})_{L^{2}(0,L)}+(\int\limits_{0}^{t}m(t-s)\mathcal{I}%
_{s}^{2}\eta (x,s)ds,\mathcal{I}_{t}\eta )_{L^{2}(0,L)}  \notag \\
&=&-\frac{d}{dt}\int\limits_{0}^{L}\mathcal{I}_{t}^{2}\eta _{x}\left(
\int\limits_{0}^{t}m\left( t-s\right) (\mathcal{I}_{s}^{2}\eta
_{x})(x,s)ds\right) dx  \notag \\
&&+\int\limits_{0}^{L}\mathcal{I}_{t}^{2}\eta _{x}\left(
\int\limits_{0}^{t}m^{\prime }\left( t-s\right) (\mathcal{I}_{s}^{2}\eta
_{x})(x,s)ds\right) dx+\int\limits_{0}^{L}m\left( 0\right) \left( \mathcal{I}%
_{t}^{2}\eta _{x}\right) ^{2}dx  \notag \\
&&+(\int\limits_{0}^{t}m(t-s)\mathcal{I}_{s}^{2}\eta (x,s)ds,\mathcal{I}%
_{t}\eta )_{L^{2}(0,L)}.  \label{e4.12*}
\end{eqnarray}%
Insertion of equations (\ref{e4.5})-(\ref{e4.12*}) into (\ref{e4.4}), and
using Lemma 2.2, yields%
\begin{eqnarray}
&&\frac{\rho _{1}}{2}^{C}\partial _{t}^{\alpha }\left\Vert \mathcal{I}%
_{t}\xi )\right\Vert _{L^{2}(0,L)}^{2}+\frac{\rho _{1}}{2}^{C}\partial
_{t}^{\alpha }\left\Vert \mathcal{I}_{x}\mathcal{I}_{t}\xi )\right\Vert
_{L^{2}(0,L)}^{2}+\frac{\kappa _{1}\partial }{2\partial t}\Vert \mathcal{I}%
_{t}^{2}\xi _{x}\Vert _{L^{2}(0,L)}^{2}  \notag \\
\ &&+\frac{\kappa _{1}\partial }{2\partial t}\Vert \mathcal{I}_{t}^{2}\xi
\Vert _{L^{2}(0,L)}^{2}+\Vert \mathcal{I}_{t}\xi \Vert _{L^{2}(0,L)}^{2}+%
\frac{\rho _{2}}{2}^{C}\partial _{t}^{\alpha }\left\Vert \mathcal{I}_{t}\eta
)\right\Vert _{L^{2}(0,L)}^{2}  \notag \\
&&+\frac{\rho _{2}}{2}^{C}\partial _{t}^{\alpha }\left\Vert \mathcal{I}_{x}%
\mathcal{I}_{t}\eta )\right\Vert _{L^{2}(0,L)}^{2}+\frac{\kappa _{2}\partial
}{2\partial t}\Vert \mathcal{I}_{t}^{2}\eta _{x}\Vert _{L^{2}(0,L)}^{2}+%
\frac{\kappa _{2}\partial }{2\partial t}\Vert \mathcal{I}_{t}^{2}\eta \Vert
_{L^{2}(0,L)}^{2}  \notag \\
&&+\int\limits_{0}^{L}h\left( 0\right) \left( \mathcal{I}_{t}^{2}\eta
_{x}\right) ^{2}dx  \notag \\
\ &\leq &\ \kappa _{1}(\mathcal{I}_{t}^{2}\eta _{x},\mathcal{I}_{t}\xi
)_{L^{2}(0,L)}-\kappa _{1}(\mathcal{I}_{t}^{2}\eta ,\mathcal{I}_{x}\mathcal{I%
}_{t}\xi )_{L^{2}(0,L)}+\Vert \mathcal{I}_{x}\mathcal{I}_{t}\xi \Vert
_{L^{2}(0,L)}^{2}  \notag \\
&&-\kappa _{1}(\mathcal{I}_{t}^{2}\xi _{x},\mathcal{I}_{t}\eta
)_{L^{2}(0,L)}+\kappa _{1}(\mathcal{I}_{t}^{2}\xi ,\mathcal{I}_{x}\mathcal{I}%
_{t}\eta )_{L^{2}(0,L)}  \notag \\
&&+\frac{d}{dt}\int\limits_{0}^{L}\mathcal{I}_{t}^{2}\eta _{x}\left(
\int\limits_{0}^{t}m\left( t-s\right) (\mathcal{I}_{s}^{2}\eta
_{x})(x,s)ds\right) dx  \notag \\
&&-\int\limits_{0}^{L}\mathcal{I}_{t}^{2}\eta _{x}\left(
\int\limits_{0}^{t}m^{\prime }\left( t-s\right) (\mathcal{I}_{s}^{2}\eta
_{x})(x,s)ds\right) dx  \notag \\
&&-(\int\limits_{0}^{t}m(t-s)\mathcal{I}_{s}^{2}\eta (x,s)ds,\mathcal{I}%
_{t}\eta )_{L^{2}(0,L)}.  \label{e4.14*}
\end{eqnarray}%
By using Cauchy $\epsilon $ -inequality, we estimate each term of the
right-hand side of previous relations to get

\begin{equation}
\kappa _{1}(\mathcal{I}_{t}^{2}\eta _{x},\mathcal{I}_{t}\xi
)_{L^{2}(0,L)}\leq \frac{\kappa _{1}}{2}\Vert \mathcal{I}_{t}^{2}\eta
_{x}\Vert _{L^{2}(0,L)}^{2}+\frac{\kappa _{1}}{2}\Vert \mathcal{I}_{t}\xi
\Vert _{L^{2}(0,L)}^{2},\   \label{e4.15}
\end{equation}%
\begin{equation}
-\kappa _{1}(\mathcal{I}_{t}^{2}\eta ,\mathcal{I}_{x}\mathcal{I}_{t}\xi
)_{L^{2}(0,L)}\leq \frac{\kappa _{1}}{2}\Vert \mathcal{I}_{t}^{2}\eta \Vert
_{L^{2}(0,L)}^{2}+\frac{\kappa _{1}}{2}\Vert \mathcal{I}_{x}\mathcal{I}%
_{t}\xi \Vert _{L^{2}(0,L)}^{2},\   \label{e4.16}
\end{equation}%
\begin{equation}
-\kappa _{1}(\mathcal{I}_{t}^{2}\xi _{x},\mathcal{I}_{t}\eta
)_{L^{2}(0,L)}\leq \frac{\kappa _{1}}{2}\Vert \mathcal{I}_{t}^{2}\xi
_{x}\Vert _{L^{2}(0,L)}^{2}+\frac{\kappa _{1}}{2}\Vert \mathcal{I}_{t}\eta
\Vert _{L^{2}(0,L)}^{2},  \label{e4.17*}
\end{equation}%
\begin{equation}
\kappa _{1}(\mathcal{I}_{t}^{2}\xi ,\mathcal{I}_{x}\mathcal{I}_{t}\eta
)_{L^{2}(0,L)}\leq \frac{\kappa _{1}}{2}\Vert \mathcal{I}_{t}^{2}\xi \Vert
_{L^{2}(0,L)}^{2}+\frac{\kappa _{1}}{2}\Vert \mathcal{I}_{x}\mathcal{I}%
_{t}\eta \Vert _{L^{2}(0,L)}^{2},  \label{e4.18*}
\end{equation}%
\begin{eqnarray}
&&-\int\limits_{0}^{L}\mathcal{I}_{t}^{2}\eta _{x}\left(
\int\limits_{0}^{t}m^{\prime }\left( t-s\right) (\mathcal{I}_{s}^{2}\eta
_{x})(x,s)ds\right) dx  \notag \\
&\leq &\frac{1}{2}\underset{0\leq t\leq T}{\sup }\left\vert m^{\prime
}\right\vert \left( 1+\frac{T^{2}}{2}\right) \Vert \mathcal{I}_{t}^{2}\eta
_{x}\Vert _{L^{2}(0,L)}^{2},  \label{e4.19*}
\end{eqnarray}%
\begin{eqnarray}
&&-(\int\limits_{0}^{t}m(t-s)\mathcal{I}_{s}^{2}\eta (x,s)ds,\mathcal{I}%
_{t}\eta )_{L^{2}(0,L)}  \notag \\
&&\frac{1}{2}\underset{0\leq t\leq T}{\sup }\left\vert m\right\vert \Vert
\mathcal{I}_{t}\eta \Vert _{L^{2}(0,L)}^{2}+\frac{T^{2}}{2}\underset{0\leq
t\leq T}{\sup }\left\vert m\right\vert \Vert \mathcal{I}_{t}^{2}\eta \Vert
_{L^{2}(0,L)}^{2}.  \label{e4.20*}
\end{eqnarray}%
By combining equality (\ref{e4.14}) and inequalities (\ref{e4.15})-(\ref%
{e4.20*}), we obtain

\begin{eqnarray}
&&^{C}\partial _{t}^{\alpha }\left\Vert \mathcal{I}_{t}\xi )\right\Vert
_{L^{2}(0,L)}^{2}+^{C}\partial _{t}^{\alpha }\left\Vert \mathcal{I}_{t}\eta
)\right\Vert _{L^{2}(0,L)}^{2}+^{C}\partial _{t}^{\alpha }\left\Vert
\mathcal{I}_{x}\mathcal{I}_{t}\xi )\right\Vert _{L^{2}(0,L)}^{2}  \notag \\
&&+^{C}\partial _{t}^{\alpha }\left\Vert \mathcal{I}_{x}\mathcal{I}_{t}\eta
)\right\Vert _{L^{2}(0,L)}^{2}+\frac{\partial }{\partial t}\Vert \mathcal{I}%
_{t}^{2}\xi _{x}\Vert _{L^{2}(0,L)}^{2}+\frac{\partial }{\partial t}\Vert
\mathcal{I}_{t}^{2}\eta _{x}\Vert _{L^{2}(0,L)}^{2}  \notag \\
\ &&+\frac{\partial }{\partial t}\Vert \mathcal{I}_{t}^{2}\xi \Vert
_{L^{2}(0,L)}^{2}+\frac{\partial }{\partial t}\Vert \mathcal{I}_{t}^{2}\eta
\Vert _{L^{2}(0,L)}^{2}+\int\limits_{0}^{L}\left( \mathcal{I}_{t}^{2}\eta
_{x}\right) ^{2}dx+\Vert \mathcal{I}_{t}\xi \Vert _{L^{2}(0,L)}^{2}  \notag
\\
&\leq &\mathcal{W}\left( \Vert \mathcal{I}_{t}\xi \Vert
_{L^{2}(0,L)}^{2}+\Vert \mathcal{I}_{t}\eta \Vert _{L^{2}(0,L)}^{2}+\Vert
\mathcal{I}_{x}\mathcal{I}_{t}\xi \Vert _{L^{2}(0,L)}^{2}+\Vert \mathcal{I}%
_{x}\mathcal{I}_{t}\eta \Vert _{L^{2}(0,L)}^{2}\right.  \notag \\
&&\left. +\Vert \mathcal{I}_{t}^{2}\xi _{x}\Vert _{L^{2}(0,L)}^{2}+\Vert
\mathcal{I}_{t}^{2}\eta _{x}\Vert _{L^{2}(0,L)}^{2}+\Vert \mathcal{I}%
_{t}^{2}\xi \Vert _{L^{2}(0,L)}^{2}+\Vert \mathcal{I}_{t}^{2}\eta \Vert
_{L^{2}(0,L)}^{2}\right) ,  \label{e4.21**}
\end{eqnarray}%
where%
\begin{equation}
\mathcal{W}=\frac{\max \left\{ \frac{\kappa _{1}}{2}+(1+\frac{T^{2}}{2})%
\frac{1}{2}\sup \left\vert m^{\prime }\right\vert ,\frac{\kappa _{1}}{2}+%
\frac{T^{2}}{2}\frac{1}{2}\sup \left\vert m\right\vert \right\} }{\min
\left\{ 1,\frac{\rho _{1}}{2},\frac{\rho _{2}}{2},\frac{\kappa _{1}}{2},%
\frac{\kappa _{2}}{2},h(0)\right\} }.  \label{e4.22*}
\end{equation}%
By discarding the last two terms from the left hand side of (\ref{e4.21*}),
replacing $t$ by $\tau $ in (\ref{e4.22*}) and then integrating with respect
to $\tau $ over the interval $(0,t),$ we obtain%
\begin{eqnarray}
&&D_{t}^{\alpha -1}\left\Vert \mathcal{I}_{t}\xi )\right\Vert
_{L^{2}(0,L)}^{2}+D_{t}^{\alpha -1}\left\Vert \mathcal{I}_{t}\eta
)\right\Vert _{L^{2}(0,L)}^{2}+D_{t}^{\alpha -1}\left\Vert \mathcal{I}_{x}%
\mathcal{I}_{t}\xi )\right\Vert _{L^{2}(0,L)}^{2}  \notag \\
&&+D_{t}^{\alpha -1}\left\Vert \mathcal{I}_{x}\mathcal{I}_{t}\eta
)\right\Vert _{L^{2}(0,L)}^{2}+\Vert \mathcal{I}_{t}^{2}\xi _{x}\Vert
_{L^{2}(0,L)}^{2}+\Vert \mathcal{I}_{t}^{2}\eta _{x}\Vert _{L^{2}(0,L)}^{2}
\notag \\
&&+\Vert \mathcal{I}_{t}^{2}\xi \Vert _{L^{2}(0,L)}^{2}+\Vert \mathcal{I}%
_{t}^{2}\eta \Vert _{L^{2}(0,L)}^{2}  \notag \\
&\leq &\mathcal{W}\left( \Vert \mathcal{I}_{t}\xi \Vert
_{L^{2}(0,t;L^{2}(0,L))}^{2}+\Vert \mathcal{I}_{t}\eta \Vert
_{L^{2}(0,t;L^{2}(0,L))}^{2}+\Vert \mathcal{I}_{x}\mathcal{I}_{t}\xi \Vert
_{L^{2}(0,t;L^{2}(0,L))}^{2}+\Vert \mathcal{I}_{x}\mathcal{I}_{t}\eta \Vert
_{L^{2}(0,t;L^{2}(0,L))}^{2}\right.  \notag \\
&&\left. +\Vert \mathcal{I}_{t}^{2}\xi _{x}\Vert
_{L^{2}(0,t;L^{2}(0,L))}^{2}+\Vert \mathcal{I}_{t}^{2}\eta _{x}\Vert
_{L^{2}(0,t;L^{2}(0,L))}^{2}+\Vert \mathcal{I}_{t}^{2}\xi \Vert
_{L^{2}(0,t;L^{2}(0,L))}^{2}+\Vert \mathcal{I}_{t}^{2}\eta \Vert
_{L^{2}(0,t;L^{2}(0,L))}^{2}\right) .  \label{e4.23*}
\end{eqnarray}%
If we omit the first four terms on the left hand side of (\ref{e4.23*}), and
use Gronwall-Bellman lemma,\ by taking%
\begin{equation}
\left\{
\begin{array}{c}
\begin{array}{c}
\mathcal{R}(t)=\Vert \mathcal{I}_{t}^{2}\xi _{x}\Vert
_{L^{2}(0,t;L^{2}(0,L))}^{2}+\Vert \mathcal{I}_{t}^{2}\eta _{x}\Vert
_{L^{2}(0,t;L^{2}(0,L))}^{2} \\
+\Vert \mathcal{I}_{t}^{2}\xi \Vert _{L^{2}(0,t;L^{2}(0,L))}^{2}+\Vert
\mathcal{I}_{t}^{2}\eta \Vert _{L^{2}(0,t;L^{2}(0,L))}^{2},%
\end{array}
\\
\frac{\partial \mathcal{R}(t)}{\partial t}=\Vert \mathcal{I}_{t}^{2}\xi
_{x}\Vert _{L^{2}(0,L)}^{2}+\Vert \mathcal{I}_{t}^{2}\eta _{x}\Vert
_{L^{2}(0,L)}^{2} \\
+\Vert \mathcal{I}_{t}^{2}\xi \Vert _{L^{2}(0,L)}^{2}+\Vert \mathcal{I}%
_{t}^{2}\eta \Vert _{L^{2}(0,L)}^{2}, \\
\mathcal{R}(t)=0,%
\end{array}%
\right.  \label{e4.24**}
\end{equation}%
we obtain%
\begin{eqnarray}
\mathcal{R}(t) &\leq &Te^{T\mathcal{W}}\left( \Vert \mathcal{I}_{t}\xi \Vert
_{L^{2}(0,t;L^{2}(0,L))}^{2}+\Vert \mathcal{I}_{t}\eta \Vert
_{L^{2}(0,t;L^{2}(0,L))}^{2}+\Vert \mathcal{I}_{x}\mathcal{I}_{t}\xi \Vert
_{L^{2}(0,t;L^{2}(0,L))}^{2}\right.  \notag \\
&&\left. +\Vert \mathcal{I}_{x}\mathcal{I}_{t}\eta \Vert
_{L^{2}(0,t;L^{2}(0,L))}^{2}\right) .  \label{e4.25*}
\end{eqnarray}%
Next, if we disregard the last four terms on the left-hand side and take
into account the inequality (\ref{e4.25*}), we end with%
\begin{eqnarray}
&&D_{t}^{\alpha -1}\left\Vert \mathcal{I}_{t}\xi )\right\Vert
_{L^{2}(0,L)}^{2}+D_{t}^{\alpha -1}\left\Vert \mathcal{I}_{t}\eta
)\right\Vert _{L^{2}(0,L)}^{2}+D_{t}^{\alpha -1}\left\Vert \mathcal{I}_{x}%
\mathcal{I}_{t}\xi )\right\Vert _{L^{2}(0,L)}^{2}  \notag \\
&&+D_{t}^{\alpha -1}\left\Vert \mathcal{I}_{x}\mathcal{I}_{t}\eta
)\right\Vert _{L^{2}(0,L)}^{2}  \notag \\
&\leq &\mathcal{W}\left( Te^{T\mathcal{W}}+1\right) \left( \Vert \mathcal{I}%
_{t}\xi \Vert _{L^{2}(0,t;L^{2}(0,L))}^{2}+\Vert \mathcal{I}_{t}\eta \Vert
_{L^{2}(0,t;L^{2}(0,L))}^{2}\right.  \notag \\
&&+\left. \Vert \mathcal{I}_{t}^{2}\xi \Vert
_{L^{2}(0,t;L^{2}(0,L))}^{2}+\Vert \mathcal{I}_{t}^{2}\eta \Vert
_{L^{2}(0,t;L^{2}(0,L))}^{2}\right) .  \label{e4.26*}
\end{eqnarray}%
Now, we are able to apply lemma 2.2, by letting%
\begin{equation}
\left\{
\begin{array}{c}
\begin{array}{c}
Q(t)=\Vert \mathcal{I}_{t}\xi \Vert _{L^{2}(0,t;L^{2}(0,L))}^{2}+\Vert
\mathcal{I}_{t}\eta \Vert _{L^{2}(0,t;L^{2}(0,L))}^{2} \\
\Vert \mathcal{I}_{t}^{2}\xi \Vert _{L^{2}(0,t;L^{2}(0,L))}^{2}+\Vert
\mathcal{I}_{t}^{2}\eta \Vert _{L^{2}(0,t;L^{2}(0,L))}^{2},%
\end{array}
\\
^{C}\partial _{t}^{\alpha }Q(t)=D_{t}^{\alpha -1}\left\Vert \mathcal{I}%
_{t}\xi )\right\Vert _{L^{2}(0,L)}^{2}+D_{t}^{\alpha -1}\left\Vert \mathcal{I%
}_{t}\eta )\right\Vert _{L^{2}(0,L)}^{2} \\
+D_{t}^{\alpha -1}\left\Vert \mathcal{I}_{x}\mathcal{I}_{t}\xi )\right\Vert
_{L^{2}(0,L)}^{2}+D_{t}^{\alpha -1}\left\Vert \mathcal{I}_{x}\mathcal{I}%
_{t}\eta )\right\Vert _{L^{2}(0,L)}^{2}, \\
Q(0)=0,%
\end{array}%
\right.  \label{e4.27*}
\end{equation}%
we infer from (\ref{e4.26*}) that%
\begin{equation}
Q(t)\leq \Gamma (\alpha )E_{\alpha ,\alpha }(\mathcal{W}\left( Te^{T\mathcal{%
W}}+1\right) t^{\alpha })D_{t}^{-\alpha }(0)=0.  \label{e4.28*}
\end{equation}%
We conclude from (\ref{e4.28*}), and (\ref{e4.27*}) that $\xi =0,$ $\eta =0.$
Consequently, $W(x,t)=(\Lambda _{1}(x,t),\Lambda _{2}(x,t))=(0,0)$ a.e in $%
Q^{T}.$

We now consider the general case for density

Since $\mathcal{E}$ is a Hilbert space, then $\overline{R(\mathcal{G})}=%
\mathcal{E}$ $\Leftrightarrow R(\mathcal{G})^{\bot }=\left\{ 0\right\}
\Leftrightarrow (\mathcal{GZ},\mathcal{K})_{\mathcal{E}}=0,$ for all $%
\mathcal{Z\in B}$ , and $\mathcal{K\in E},$then $\mathcal{K=(K}_{1},\mathcal{%
K}_{2})=\left\{ (J_{1},J_{3},J_{4}),(J_{2},J_{5},J_{6})\right\} =(0,0),$
that is $J_{1}=J_{2}=J_{3}=J_{4}=J_{5}=J_{6}=0.$ So suppose that for some
element $\mathcal{K=(K}_{1},\mathcal{K}_{2})=\left\{
(J_{1},J_{3},J_{4}),(J_{2},J_{5},J_{6})\right\} \in R(\mathcal{G})^{\bot }$%
\begin{eqnarray}
&&(\mathcal{GZ},\mathcal{K})_{\mathcal{E}}  \notag \\
&=&(\{\mathcal{S}_{1}(\mathcal{\theta },\phi ),\mathcal{S}_{2}(\mathcal{%
\theta },\phi \},\{\mathcal{K}_{1},\mathcal{K}_{2}\})_{\mathcal{E}}  \notag
\\
&=&(\{\mathcal{S}_{1}(\mathcal{\theta },\phi ),\Gamma _{1}\mathcal{\theta }%
,\Gamma _{2}\mathcal{\theta }\},\{\mathcal{S}_{2}(\mathcal{\theta },\phi
),\Gamma _{1}\phi ,\Gamma _{2}\phi
\}\},\{\{J_{1},J_{2},J_{3}\},\{J_{4},J_{5},J_{6}\}\})_{\mathcal{E}}  \notag
\\
&=&(\mathcal{S}_{1}(\mathcal{\theta },\phi ),J_{1})_{L^{2}(Q^{T})}+(\Gamma
_{1}\mathcal{\theta },J_{2})_{L^{2}(0,L)}+(\Gamma _{2}\mathcal{\theta }%
,J_{3})_{L^{2}(0,L)}  \notag \\
&&+(\mathcal{S}_{2}(\mathcal{\theta },\phi ),J_{4})_{L^{2}(Q^{T})}+(\Gamma
_{1}\phi ,J_{5})_{L^{2}(0,L)}+(\Gamma _{2}\phi ,J_{6})_{L^{2}(0,L)}  \notag
\\
&=&0,  \label{e4.31}
\end{eqnarray}%
where $\mathcal{Z}$ runs over the space $\mathcal{B}$, we have to prove that
$\mathcal{K}=0.$

Let $\mathcal{Z}\in D_{0}(\mathcal{G}),$ then equation (\ref{e4.31}) becomes%
\begin{equation}
(\mathcal{S}_{1}(\mathcal{\theta },\phi ),J_{1})_{L^{2}(Q^{T})}+(\mathcal{S}%
_{2}(\mathcal{\theta },\phi ),J_{4})_{L^{2}(Q^{T})}=0.  \label{e4.32}
\end{equation}%
Hence, by virtue of Theorem 4.2, it follows from (\ref{e4.32}) that $%
J_{1}=J_{4}=0$. Consequently, equation (\ref{e4.31}) takes the form%
\begin{equation}
(\Gamma _{1}\mathcal{\theta },J_{2})_{L^{2}(0,L)}+(\Gamma _{2}\mathcal{%
\theta },J_{3})_{L^{2}(0,L)}+(\Gamma _{1}\phi ,J_{5})_{L^{2}(0,L)}+(\Gamma
_{2}\phi ,J_{6})_{L^{2}(0,L)}=0.  \label{e4.33}
\end{equation}%
Since the four terms in (\ref{e4.33}) vanish independently and since the
ranges $R(\Gamma _{1}),R(\Gamma _{2})$ of the trace operators $\Gamma
_{1},\Gamma _{2}$ are everywhere dense in the space $L^{2}(0,L)$, then it
follows from (\ref{e4.33}) that $J_{2}=J_{3}=J_{5}=J_{6}=0$. Consequently $%
\mathcal{K}=0$, that is $R(\mathcal{G})^{\perp }=\{0\}$. Thus $\overline{R(%
\mathcal{G})}=\mathcal{E}$.

\bigskip

\textbf{Conclusion}: In this article, we proved the well posedness of a
nonhomogeneous Timoshenko system with a viscoelastic damping term.The
coupled two hyperbolic equations were associated with initial conditions and
a nonlocal boundary conditions. The proofs of the results are mainly based
on some energy and a priori estimates and on some density arguments. The
method uses functional analysis tools such as operator theory and density
arguments. It is found that the method is efficient and powerful for solving
initial boundary value problems with non-local constraints.The a priori
estimate for the solution can be provided by constructing a suitable
multiplicator and from which it is also possible to establish the
solvability of the stated problem. We notice that no previous works were
treated for Timoshenko systems with nonlocal conditions of integral type.

\textbf{Acknowledgment}: The authors would like to extend their sincere
appreciation to the Deanship of Scientific Research at King Saud University
for its funding this Research group No (RG 117).

The second author would like to thank Shagrah University for giving her the
chance to further her PhD studies at King Saud University.

\textbf{Availability of data and materials: \ }Data sharing is not
applicable to this article as no data sets were generated or analysed during
the current study.

\textbf{Competing interests: }The authors declare that they have no
competing interests.

\textbf{Authors' contributions: \ }All authors contributed equally to this
work. All authors read and approved the final manuscript.

\end{document}